\documentclass[12pt]{amsart}
\usepackage{amsthm,amsmath,amssymb}
\usepackage{dsfont}
\usepackage{mathrsfs} 
\usepackage[all,cmtip]{xy}
\usepackage{patchcmd}
\input{diagxy}
\usepackage{url}
\usepackage{bm}
\usepackage{mathtools}
\usepackage{enumerate}
\usepackage{pbox}
\usepackage{color}
\usepackage{stmaryrd}

 \theoremstyle{plain}
 \newtheorem{thm}{Theorem}[section]
 \newtheorem{cor}[thm]{Corollary}
 \newtheorem{lem}[thm]{Lemma}
 
 \newtheorem{prop}[thm]{Proposition}

\theoremstyle{definition}
 \newtheorem{defn}[thm]{Definition}

\theoremstyle{remark}
 \newtheorem{rem}[thm]{Remark}

 \newtheorem{nota}[thm]{Notation}
 \newtheorem{conv}[thm]{Convention}
 \newtheorem{exam}[thm]{Example}

 \numberwithin{equation}{section}

\newtheorem*{ack}{Acknowledgment}

\DeclareMathOperator{\VF}{VF}
\DeclareMathOperator{\ACVF}{ACVF}
\DeclareMathOperator{\RV}{RV}
\DeclareMathOperator{\DC}{DC}
\DeclareMathOperator{\MM}{\mathcal{M}}

\DeclareMathOperator{\OO}{\mathcal{O}}

\DeclareMathOperator{\UU}{\mathcal{U}}

 \DeclareMathOperator{\dom}{dom}

 \DeclareMathOperator{\id}{id}

 \DeclareMathOperator{\lh}{lh}

 \DeclareMathOperator{\cha}{char}
 
 \DeclareMathOperator{\supp}{supp}

 \DeclareMathOperator{\acl}{acl}
 
 \DeclareMathOperator{\pr}{pr}

 \DeclareMathOperator{\mgl}{GL}
\DeclareMathOperator{\msl}{SL}
\DeclareMathOperator{\jcb}{Jcb}

\DeclareMathOperator{\K}{\Bbbk}

\def\Xint#1{\mathchoice
{\XXint\displaystyle\textstyle{#1}}%
{\XXint\textstyle\scriptstyle{#1}}%
{\XXint\scriptstyle\scriptscriptstyle{#1}}%
{\XXint\scriptscriptstyle\scriptscriptstyle{#1}}%
\!\int}
\def\XXint#1#2#3{{\setbox0=\hbox{$#1{#2#3}{\int}$}
\vcenter{\hbox{$#2#3$}}\kern-.5\wd0}}

\newcommand{\Q}{\mathds{Q}}
\newcommand{\N}{\mathds{N}}

\newcommand{\R}{\mathds{R}}

\newcommand{\p}{$p$\nobreakdash}
\newcommand{\omin}{$o$\nobreakdash}

\newcommand{\gC}{\mathfrak{C}}
\newcommand{\gD}{\mathfrak{D}}

\newcommand{\ga}{\mathfrak{a}}
\newcommand{\gb}{\mathfrak{b}}
\newcommand{\gc}{\mathfrak{c}}
\newcommand{\gd}{\mathfrak{d}}

\newcommand{\gf}{\mathfrak{f}}

\newcommand{\go}{\mathfrak{o}}
\newcommand{\gp}{\mathfrak{p}}

\newcommand{\0}{\emptyset}

\DeclareMathAlphabet{\mathpzc}{OT1}{pzc}{m}{it}

 \newcommand{\abs}[1]{\left\vert#1\right\vert}
 
 \newcommand{\set}[1]{\left\{#1\right\}}

 \newcommand{\wh}[1]{\widehat{#1}}
 
 \newcommand{\dhat}[1]{\wh{\wh{#1}}}

\newcommand{\mdl}[1]{\mathcal{#1}}  
\newcommand{\bb}[1]{\mathbb{#1}}

\newcommand{\ex}[1]{\exists #1 \;} 

\newcommand{\rest}{\upharpoonright}

\newcommand{\fun}{\longrightarrow}
\newcommand{\efun}{\longmapsto}
\newcommand{\sub}{\subseteq}

\newcommand{\mi}{\smallsetminus}

\DeclareMathOperator{\mVF}{\mu \! \VF}

\DeclareMathOperator{\mRV}{\mu \! \RV}

\DeclareMathOperator{\mG}{\mu \Gamma}

\DeclareMathOperator{\mRES}{\mu RES}

\DeclareMathOperator{\rv}{rv}

\DeclareMathOperator{\csn}{csn}
\DeclareMathOperator{\rcsn}{\overline {csn}}
\DeclareMathOperator{\vv}{val}

\DeclareMathOperator{\gsk}{\mathbf{K}_+}
\DeclareMathOperator{\ggk}{\mathbf{K}}

\DeclareMathOperator{\ob}{Ob}
\DeclareMathOperator{\fn}{FN}
\DeclareMathOperator{\fib}{fib}
\DeclareMathOperator{\fin}{fin}
\DeclareMathOperator{\vol}{vol}
\DeclareMathOperator{\isp}{I_{sp}}
\DeclareMathOperator{\misp}{\mu I_{sp}}

\DeclareMathOperator{\rad}{rad}

\DeclareMathOperator{\vrv}{vrv}

\DeclareMathOperator{\RVH}{RVH}

\DeclareMathOperator{\KRC}{\hat{\mathbf{K}} \mathds R}
\DeclareMathOperator{\KCC}{\hat{\mathbf{K}} \mathds C}

\DeclareMathOperator{\tbk}{tbk}

\DeclareMathOperator{\res}{res}

\DeclareMathOperator{\leb}{\Lambda}

\DeclareMathOperator{\ifn}{Int}

\newcommand{\lan}[1]{\mathcal{L}_{\textup{#1}}}
\DeclareMathOperator{\KRCV}{\mu \! \RV}
\DeclareMathOperator{\KRCES}{\mu RES}

\author[Yimu Yin]{Yimu Yin}
\address{Institut Math\'{e}matique de Jussieu \\ Universit\'e Pierre et Marie Curie \\ 4 place Jussieu \\ 75252 Paris Cedex 05 \\ France}
\email{yyin@math.jussieu.fr}
\title[Fourier transform in $\ACVF$]{Fourier transform of the additive group in algebraically closed valued fields}

\begin{document}

\begin{abstract}
We continue the study of the Hrushovski-Kazhdan integration theory and consider exponential integrals. The Grothendieck ring is enlarged via a tautological additive character and hence can receive such integrals. We then define the Fourier transform in our integration theory and establish some fundamental properties of it. Thereafter a basic theory of distributions is also developed. We construct the Weil representations in the end as an application. The results are completely parallel to the classical ones.
\end{abstract}

\subjclass{Primary 03C60; Secondary 11S80, 20C08}
\keywords{Fourier transform, motivic integration, algebraically closed valued fields}

\maketitle

\tableofcontents

\section{Introduction}

The Hrushovski-Kazhdan integration theory~\cite{hrushovski:kazhdan:integration:vf} is a major development in the theory of motivic integration. The fundamental idea of the theory is to construct canonical homomorphisms between various Grothendieck rings associated with the first-order theory $\ACVF$ of algebraically closed valued fields. The simplest of these constructions are presented in~\cite{Yin:special:trans, Yin:int:acvf}. In this paper, based on the work and ideas in~\cite{hrushovski:kazhdan:integration:vf}, we demonstrate how to extend the construction to include parametrized exponential integrals, typically of the form
\[
\int_{y \in \VF^m} \int_{x \in \VF^n} \bm f( x,  y) \exp(\bm g( x, y)),
\]
where the requirements on the definable functions $\bm f$, $\bm g$ are very natural.

To describe in a few words how motivic integration is different from classical integration it seems best to begin by pointing out that the ring that provides values for integrals is not the real field but a Grothendieck ring. The latter is traditionally constructed from equivalence classes of algebraic varieties and, more generally in the model-theoretic setting, from equivalence classes of definable subsets. Topological tools that are essential to many classical constructions are no longer available; instead, since it was first introduced by Maxim Kontsevich in 1995, techniques from first-order model theory of definable sets underlie much of the development of this new kind of integration. In fact, at risk of being overly simple-minded, one may think of motivic integration as classical integration with the topological concepts of ``continuity'', ``convergence'', etc.\ replaced everywhere by the model-theoretic concept of ``definability''.

To be sure, the class of definable integrals is conceptually narrower than the class of integrals that can be more or less dealt with classically. However, there are many reasons why the motivic approach to integration will play an increasingly important role. We mention two here.

Firstly, the progress in model theory in the last few decades suggests that many natural mathematical properties are subject to first-order treatment. In our context, given the fact that some very complicated integral identities are already motivic (see, for example, \cite{cluckers:hales:loeser:transfer, cunning:hales:good}), it is reasonable to expect that many other important kinds of integrals are definable in some first-order languages and hence may be studied motivically. We note that, in their recent paper~\cite{hru:kazh:2009}, Hrushovski and Kazhdan have developed a partially first-order method to study adelic structures over curves and, in particular, have obtained a global Poisson summation formula.

Secondly, if one is more interested in the structure of a space of functions (for example, functional equations) than actual computation of functions, then constantly worrying about things such as convergence seems to be an unnecessary burden. By this we just mean that there is no need to insist on assigning ``numerical values'' to integrals, especially when it is not possible, and sometimes working with ``geometric values'' is more effective. Definable integrals are of a more geometric nature and are better behaved, at least before specializing to local fields. Some pathological phenomena afforded by point-set topology are thus avoided. For example, while classically it is possible that two iterative integrals of a function exist but are not equal, this cannot happen to definable integrals. This is our version of the Fubini theorem (Proposition~\ref{fubini}).

This better behavior of definable integrals mentioned above may be a result of how motivic measure is manufactured: the volume of a geometric object is somehow provided by the object itself, subject to certain geometric equivalence relations. The construction of exponential integrals in this paper is very illustrative of this ``tautological'' nature. To expand on this point, let us take a brief moment to outline how canonical homomorphisms between various Grothendieck rings are constructed in \cite{hrushovski:kazhdan:integration:vf,Yin:special:trans, Yin:int:acvf,Yin:int:expan:acvf}.

Let $(K, \vv : K \fun \Gamma)$ be an algebraically closed valued field, where $\vv$ is the valuation map, and $\OO$, $\MM$, $\K$ the corresponding valuation
ring, its maximal ideal, and the residue field. Let
\[
\RV(K) = K^{\times} / (1 + \MM)
\]
and $\rv : K^{\times} \fun \RV(K)$ be the quotient map. Note that, for each $a \in K$, $\vv$ is constant on the subset $a + a\MM$ and hence there is a naturally induced map $\vrv$ from $\RV(K)$ onto the value group $\Gamma$. The situation is illustrated in the following commutative diagram
\begin{equation*}
\bfig
 \square(0,0)/^{ (}->`->>`->>`^{ (}->/<600, 400>[\OO \mi \MM`K^{\times}`\K^{\times}`
\RV(K);`\text{quotient}`\rv`]
 \morphism(600,0)/->>/<600,0>[\RV(K)`\Gamma;\vrv]
 \morphism(600,400)/->>/<600,-400>[K^{\times}`\Gamma;\vv]
\efig
\end{equation*}
where the bottom sequence is exact. This structure may be expressed by a two-sorted first-order language $\lan{RV}$ (see Defintion~\ref{defn:lrv}), where $K$ is referred to as the $\VF$-sort and $\RV$ is taken as a new sort, the $\RV$-sort. Moreover, there could be extra structure in the diagram above; in particular, for the construction in this paper, there is a cross-section $\csn : \Gamma \fun K^{\times}$ (see Defintion~\ref{defn:cross}).

Now let $\mVF[*]$ and $\mRV[*]$ be two categories of definable sets with volume forms that are respectively associated with the $\VF$-sort and the $\RV$-sort, where the notation ``$[*]$'' means gradation by ambient dimension. The main construction of the Hrushovski-Kazhdan theory is a canonical homomorphism from the Grothendieck semiring $\gsk \mVF[*]$ to the Grothendieck semiring $\gsk \mRV[*]$ modulo a semiring congruence relation $\isp$ on the latter. In fact, it turns out to be an isomorphism. This construction has three main steps.
\begin{enumerate}[{Step} 1.]
 \item First we define a lifting map $\bb L$ from the set of objects in $\mRV[*]$ into the set of objects in $\mVF[*]$. Next we single out a subclass of isomorphisms in $\mVF[*]$, which are called special bijections. Then we show that for any object $A$ in $\mVF[*]$ there is a special bijection
$T$ on $A$ and an object $\bm U$ in $\mRV[*]$ such that $T(A)$ is isomorphic to $\bb L (\bm U)$. This implies that $\bb L$
hits every isomorphism class of $\mVF[*]$. Of course, for this result alone we do not have to limit our means to special
bijections. However, in Step~3 below, special bijections become an essential ingredient in computing the semiring congruence
relation $\isp$.

 \item For any two isomorphic objects $\bm U_1$, $\bm U_2$ in $\mRV[*]$, their lifts $\bb L(\bm U_1), \bb L(\bm U_2)$ in
$\mVF[*]$ are isomorphic as well. This shows that $\bb L$ induces a semiring homomorphism from $\gsk \mRV[*]$ into $\gsk \mVF[*]$, which is also denoted by $\bb L$.

 \item A number of classical properties of integration can already be verified for the inversion of the homomorphism $\bb L$ and hence, morally, this third step is not necessary. For applications, however, it is much more satisfying to have a precise description of the semiring congruence relation induced by $\bb L$. The basic notion used in the description is that of a blowup (or dilatation) of an object in $\mRV[*]$, which is essentially a restatement of the trivial fact that there is an additive translation from $1 + \MM$ onto $\MM$. We then show that, for any objects $\bm U_1$, $\bm U_2$ in $\mRV[*]$, there are isomorphic blowups $\bm U_1^{\sharp}$, $\bm U_2^{\sharp}$ of them if and only if $\bb L(\bm U_1)$, $\bb L(\bm U_2)$ are isomorphic. The ``if'' direction essentially contains a form of Fubini's Theorem and is the most technically involved part of the construction.
\end{enumerate}
The inverse of $\bb L$ thus obtained is called a Grothendieck semiring homomorphism. If the Jacobian transformation preserves integrals, that is, the change of variables formula holds, then it may be called a motivic integration; when the semirings are formally groupified, it is recast as a ring homomorphism, which is denoted by $\int$ and its target ring by $\KRC$. See Theorem~\ref{main:prop:k:vol:dag} for the particular version that we shall use in this paper.

The integration formalism just described is unable to accommodate additive characters because $\KRC$ is not big enough. To remedy this, we just take a quotient $\Omega$ of the additive structure of the field $\VF$ and add it to $\KRC$, as a set of symbols, to form a group ring $\KCC$ (hence additive relation is turned into multiplicative relation).

More precisely, let $\Omega = \VF / \MM$ and $\theta : \VF \fun \Omega$ be the quotient map. There is a natural map $\Omega \fun \Gamma$ induced by $\vv$, which will also be denoted by $\vv$. For any \p-adic field $\Q_p$, the specialization $\Omega(\Q_p)$ of $\Omega$ to $\Q_p$ may be identified naturally with the $p$th power roots of unity via any additive character $\chi$ such that, for any \p-adic number $a$ with $\vv(a) \leq 0$, $\chi(a)$ is a $p^{- \vv(a) + 1}$th root of unity. In general, let $F$ be a number field, $\vv_p$ a valuation of $F$ with respect to a prime number $p$, $F_p$ the completion of $F$ with respect to $\vv_p$, and $\OO_p$, $\MM_p$ the corresponding valuation ring and its maximal ideal. Let $\chi$ be an additive character of $F_p$. The \emph{conductor} $\gf_{\chi}$ of $\chi$ is the largest ideal of the form $\MM_p^m$, where $m$ is a nonnegative integer ($\MM_p^0 = \OO_p$ by definition), such that $\chi(\MM_p^m) = 1$. This always exists, see \cite[Lemma~4, p.~114]{weil:basic:number}. The integer $m$ is the \emph{multiplicity} of $\gf_{\chi}$. Now, if we view $\Omega$ as a motivic analogue of the subgroup of roots of unity of the unit circle and $\theta$ a motivic analogue of a generic additive character, then motivic integration with respect to $\theta$, when specialized to (non-archimedean) completions of $F$, may be viewed as integration with respect to all additive characters with multiplicity 1 for almost all $p$ at once.

There are certain fundamental properties of additive characters that must hold in this group ring $\KCC$. To achieve these we can take quotients or localize or employ other standard algebraic operations. As one can easily surmise, this sort of construction is very flexible and many other desirable features may be incorporated along similar lines.

Let us now describe in more detail how the sections of this paper are organized. In Section~\ref{section:1st:prop} we introduce and recall notation and terminology. There are also a few technical results that will be needed in later sections. They may be skipped (at least the proofs) without impairing the essential understanding of the rest of the paper and hence are deferred to the last section. Some simple modifications of the target ring of integration are also carried out. This ring, in this paper, is our analogue of the real field. It is no longer graded and has more invertible elements. It is further enlarged in Section~\ref{section:integrable} via a tautological additive character and thereby becomes an analogue of the complex field. We then single out a subclass of definable functions that will be the focus of the discussions in the subsequent sections, namely integrable functions. There is flexibility in the concept of integrability and the class can be made larger. But the cutoff line we have adopted seems most natural: desirable features such as the Fubini property and closure under convolution can be proved easily. In Section~\ref{section:with:add:char} the Fourier transform is defined, which perhaps should be called the Fourier-Laplace transform. Thereafter various fundamental integral identities involving the Fourier transform are established, for example, the convolution formula, the Fourier inversion formula, the Plancherel formula, etc. Although distributions may seem manifestly non-first-order, a basic theory of definable distributions can be developed within our framework, which we shall do in Section~\ref{section:dist}. In the last section, as an application, we show that the Weil representations exist on the Schwartz spaces associated with algebraically closed valued fields. For simplicity, we shall only show this for $\msl_2$. It is not hard to see that in principle the construction goes through for all symplectic groups, but the computations involved are too complicated to be presented in a clear manner.

We remark that the generalization of the construction of Weil representations in this paper is a na\"ive one. It merely creates a formalism in which the classical arguments make sense in the context of algebraically closed valued fields, although some interesting phenomena have already appeared along the way. Ideally the generalized Weil representations over algebraically closed valued fields should be the ``limit'' of the classical Weil representations over locally compact fields. To make this connection we need to incorporate data from the Galois group. This will be explored in a sequel.

On the level of local fields there is also a very general approach to motivic integration, namely the Cluckers-Loeser theory~\cite{cluckers:loeser:constructible:motivic:functions}; see~\cite{gordon:yaffe:2008} for an excellent exposition. Their construction of exponential integrals is contained in~\cite{cluckers:loeser:motivic:fourier:transform}. There is also a comparable theory of motivic integration for real closed fields (see \cite{Yin:int:tcvf}). For a general introduction to the development before these general approaches we refer to the articles~\cite{hales:2005, loo:mot}.

\begin{ack}
I would like to thank Thomas Hales for many hours of stimulating discussions during the preparation of this paper, in particular, for pointing out to me the possibility of constructing a motivic version of the Weil representations. The research reported in this paper has been partially supported by the ERC Advanced Grant NMNAG.
\end{ack}

\section{Preliminaries}\label{section:1st:prop}

The construction of exponential integrals in this paper is of a formal nature and works for any version of the Hrushovski-Kazhdan style motivic integration that has been developed so far: with or without volume forms, with or without sections, with or without arbitrary extra structure in the $\RV$-sort, and so forth. For ease of discussion and concreteness, we shall mainly concentrate on a scenario that is treated in \cite{Yin:int:expan:acvf}, namely the one that extends an algebraically closed valued field, considered as a model of the $\lan{RV}$-theory $\ACVF$, by a cross-section and thereby introduces volume forms on objects over the residue field.

The reader is referred to~\cite{Yin:QE:ACVF:min, Yin:int:acvf, Yin:special:trans, Yin:int:expan:acvf} for notation and terminology. There will be reminders and, inevitably, repetitions as we go along. To begin with, the nomenclature and the notation in \cite[Definitions~2.1, 2.2, 2.6, Notation~2.9]{Yin:special:trans} shall be used frequently, as well as the various notational conventions concerning coordinate projection maps in~\cite[Notation~2.10]{Yin:special:trans}.

\begin{defn}\label{defn:lrv}
The language $\lan{RV}$ has the following sorts and symbols:
\begin{itemize}
 \item a $\VF$-sort, which uses the language of rings
 $\lan{R} = \set{0, 1, +, -, \times}$;
 \item an $\RV$-sort, which uses
  \begin{itemize}
    \item the group language $\set{1, \times}$,
    \item a constant symbol $\infty$,
    \item a unary predicate $\K^{\times}$,
    \item a binary function $+ : \K^2 \fun \K$ and a unary function $-: \K \fun \K$, where $\K = \K^{\times} \cup \{\infty\}$,
    \item a binary relation $\leq$;
    \end{itemize}
  \item a function symbol $\rv$ from the $\VF$-sort into the $\RV$-sort.
\end{itemize}
We write $\VF \mi \{0\}$ and $\RV \mi \{\infty\}$ as $\VF^{\times}$ and $\RV^{\times}$, respectively.
\end{defn}

\begin{defn}\label{defn:acvf}
\emph{The theory $\ACVF$ of algebraically closed valued fields in $\lan{RV}$} states the following:
\begin{itemize}
 \item $(\VF, 0, 1, + , -, \times)$ is an algebraically close field;

 \item $(\RV^{\times}, 1, \times)$ is a divisible abelian
 group, where multiplication $\times$ is augmented by $t \times \infty =
 \infty$ for all $t \in \RV$;

 \item $(\K, \infty, 1, +, -, \times)$ is an algebraically closed field (note that the symbol $\infty$ serves as the element $0$ in $\K$, which, for psychological reasons, shall indeed be written as $0$ when $\K$ is concerned);

 \item the relation $\leq$ is a preordering on $\RV$ with
 $\infty$ the top element and $\K^{\times}$ the equivalence
class of 1;

 \item the quotient $\RV / \K^{\times}$, denoted by $\Gamma \cup \set{\infty}$, is a divisible ordered abelian
group with a top element, where the ordering and the group
operation are induced by $\leq$ and $\times$, respectively,
and the quotient map $\RV \fun \Gamma_{\infty} \coloneqq \Gamma \cup \set{\infty}$ is
denoted by $\vrv$;

 \item the function $\rv : \VF^{\times} \fun \RV^{\times}$
 is a surjective group homomorphism augmented by $\rv(0) =
\infty$ such that the composite function
\[
\vv = \vrv \circ \rv : \VF \fun \Gamma_{\infty}
\]
is a valuation with the valuation ring $\OO = \rv^{-1}(\RV^{\circ})$ and its maximal ideal $\MM =
\rv^{-1}(\RV^{\circ\circ})$, where
\[
\RV^{\circ} = \set{x \in \RV: 1 \leq x}, \quad
\RV^{\circ\circ} = \set{x \in \RV: 1 < x}.
\]
\end{itemize}
\end{defn}

\begin{defn}\label{defn:cross}
A \emph{cross-section} of $\Gamma$ is a group homomorphism $\csn : \Gamma \fun \VF^{\times}$ such that $\vv \circ \csn = \id$. The \emph{reduced cross-section} of $\Gamma$ is the function $\rcsn = \rv \circ \csn : \Gamma \fun \RV^{\times}$. Set $\csn(\infty) = 0$ and $\rcsn(\infty) = \infty$. We may and shall treat $\csn$ and $\rcsn$ as functions on $\RV$ in the natural way. The \emph{twistback} function $\tbk : \RV \fun \K$ is given by $u \efun u / \rcsn(u)$, where we set $\infty / \infty = 0$.
\end{defn}

The expansion of $\lan{RV}$ with the function symbol $\csn$ is denoted by $\lan{RV}^{\csn}$. The theory $\ACVF^{\csn}$ in $\lan{RV}^{\csn}$ states that, in addition to the axioms of $\ACVF$, $\csn$ is a cross-section of $\Gamma$. If the characteristics are specified then we write $\ACVF^{\csn}(0, p)$, etc. The theory $\ACVF^{\csn}(0, 0)$ is complete and admits quantifier elimination (see \cite[Theorem~2.6]{Yin:int:expan:acvf}). From this moment forth, we fix a sufficiently saturated model $\gC^{\csn} \models \ACVF^{\csn}(0, 0)$. Note that $\gC^{\csn}$ is a specialization of the model $\gC^{\tilde{\bb T}}$ considered in \cite{Yin:int:expan:acvf} (see the discussion after \cite[Theorem~2.6]{Yin:int:expan:acvf}). The $\lan{RV}$-reduct of $\gC^{\csn}$ is denoted by $\gC$.

\begin{nota}\label{indexing}
Unless otherwise specified, by writing $a \in A$ we shall mean that $a$ is a finite tuple $(a_1, a_2, \ldots)$ of elements (or ``points'') of $A$, whose length, denoted by $\lh(a)$, is not always indicated.

We say that $B$ is a subset \emph{in} $A$ if $B \sub A^n$ for some $n$.

For each $n \in \N$, let $[n] = \set{1, \ldots, n}$. Let $A$ be a subset of $\VF^n \times \RV^m$. As a general rule, the coordinates of $A$ are indexed by $[n+m] = [n] \uplus [m]$. Let $E \sub [n+m]$ and $\tilde E = [n+m] \mi E$. If $E$ is a singleton $\{i\}$ then we always write $E$ as $i$ and $\tilde E$ as $\tilde i$. The projection of $A$ to the coordinates in $E$ is denoted by $\pr_E(A)$. For $a \in \pr_{\tilde E} (A)$, the fiber $\{b : ( b, a) \in A \}$ is denoted by $\fib(A, a)$ or, if there is no danger of confusion, by $A_a$. Note that the distinction between the two subsets $A_a$ and $A_a \times \{ a \}$ is often immaterial and hence they will be tacitly identified. Also, it is more convenient to use simple descriptions as subscripts. For example, if $E = [k]$, etc., then we may write $\pr_{\leq k}$, etc. If $\pr_{E}(A) \sub B$ and $B$ has been clearly understood in the context then it is more informative to write $\pr_B \rest A$.

Given a function $f : A \fun B$, we shall often write $A_b$ for the pullback fiber over $b \in B$ under $f$. In particular, given a subset $A$, we may write $A_{x}$ for the fiber over $x$ under a function of the forms $\rv \rest A$, $\vv \rest A$, $\vrv \rest A$, etc. Of course which function is being considered should always be clear in context.
\end{nota}

\begin{conv}\label{conv:imag}
Let $\textup{res} : \RV \fun \K$ be the function given by $\textup{res} \rest \K^{\times} = \id$ and $\textup{res}(t) = 0$ for all $t \notin \K^{\times}$. Technically speaking, $+ : \K^2 \fun \K$ is a function symbol only in the imaginary sort $\K$, which, as in \cite{hrushovski:kazhdan:integration:vf, Yin:special:trans, Yin:int:acvf, Yin:int:expan:acvf}, is subsumed into the $\RV$-sort. An \emph{imaginary $\K$-term} is a term of the form
\[
\sum_{i = 1}^k \res(\rv(F_{i}(X)) \cdot r_{i} \cdot  Y^{ n_i}),
\]
where $X$ are $\VF$-sort variables, $Y$ are $\RV$-sort variables, $n_i \in \N$, $r_i \in \RV$, and $F_{i}(X)$ is a polynomial with coefficients in $\VF$. An \emph{imaginary $\Gamma$-term} is a term of the same form with $\res$ replaced by $\vrv$. At times it will be more convenient to treat these terms (and other similar ones) as real terms in the language $\mdl L_{\K \Gamma}$ (respectively $\mdl L_{\K \Gamma}^{\rcsn}$) that naturally corresponds to the three-sorted structure of the reduct of $\gC$ (respectively $\gC^{\csn}$) to the $\RV$-sort.
\end{conv}

\begin{defn}\label{defn:disc}
A subset $\gb$ of $\VF$ is an \emph{open disc} if there is a $\gamma \in \Gamma$ and a $b \in \gb$ such that $a \in \gb$ if and
only if $\vv(a - b) > \gamma$; it is a \emph{closed disc} if $a \in \gb$ if and only if $\vv(a - b) \geq \gamma$; it is an
\emph{$\RV$-disc} if $\gb = \rv^{-1}(t)$ for some $t \in \RV$. The value $\gamma$ is the \emph{valuative radius} of $\gb$, which is denoted by $\rad (\gb)$. Each point in $\VF$ is a closed disc of valuative radius $\infty$ and $\VF$ is a clopen disc of radius $- \infty$.

A closed disc with a maximal open subdisc removed is called a \emph{thin annulus}.

A subset $\gp \sub \VF^n \times \RV^m$ is an (\emph{open, closed, $\RV$-}) \emph{polydisc} if it is of the form $(\prod_{i \leq n} \gb_i) \times \{t \}$, where each $\gb_i$ is an (open, closed, $\RV$-) disc. The \emph{radii} of $\gp$, denoted by $\rad(\gp)$, means the tuple $((\rad(\gb_1), \ldots, (\rad(\gb_n))$, while the \emph{radius} of $\gp$ means $\min \rad(\gp)$. The open and closed polydiscs centered at $a = (a_1, \ldots, a_n) \in \VF^n$ with radii $\gamma = (\gamma_1, \ldots, \gamma_n) \in \Gamma^n$ are denoted by $\go(a, \gamma)$ and $\gc(a, \gamma)$, respectively.

An $\RV$-polydisc $\rv^{-1}(t_1, \ldots, t_n) \times \{ s \}$ is \emph{degenerate} if $t_i = \infty$ for some $i$.

The \emph{$\RV$-hull} of a subset $A$, denoted by $\RVH(A)$, is the union of all the $\RV$-polydiscs whose intersections with $A$ are nonempty. If $A$ equals $\RVH(A)$ then $A$ is called an \emph{$\RV$-pullback}.
\end{defn}

\begin{conv}\label{conv:s}
For convenience and without loss of generality, by a substructure we shall always mean a substructure that is equal to its definable closure. Let $S$ be a small substructure of $\gC^{\csn}$. Note that the $\lan{RV}$-reduct of $S$ is $\VF$-generated, which, for simplicity, shall also be denoted by $S$ if there is no danger of confusion. The corresponding expanded languages (with constants in $S$) are still referred to as $\lan{RV}$ and $\lan{RV}^{\csn}$. Parameters from $S$ are allowed and they will not be specified unless it is necessary. So in effect we shall be working with the complete theories $\ACVF(S)$ and $\ACVF^{\csn}(S)$. By an $\lan{RV}$-definable (respectively $\lan{RV}^{\csn}$-definable, etc.) subset we mean an $S$-$\lan{RV}$-definable (resp.\ $S$-$\lan{RV}^{\csn}$-definable, etc.) subset. In general, \emph{by a definable subset we mean an $\lan{RV}^{\csn}$-definable subset}, unless indicated otherwise in context. Parameters from sources other than $S$ will be specified in context.

At times it will be convenient to work in the traditional expansion $\gC^{\csn, eq}$ of $\gC^{\csn}$ when imaginary elements are called for. However, a much simpler expansion $\gC^{\csn, \bullet}$ suffices: it has only one additional sort $\DC$ that contains, as elements, all the open and closed discs. This means that, when we work in $\gC^{\csn, \bullet}$, the underlying substructure $S$ may contain discs with valuative radii in $\Gamma(S)$. These discs can be used as parameters as well. Note that it is redundant to include in $\DC$ discs centered at $0$, since they may be identified with their valuative radii. For convenience, we do think of $\VF$ as a subset of $\DC$. This expansion can help reduce the technical complexity of our discussion. However, as is the case with $\Gamma$, it is conceptually inessential since, for the purposes of this paper, all allusions to discs as (imaginary) elements may be eliminated in favor of objects already definable in $\gC^{\csn}$.

For a disc $\ga \sub \VF$ the corresponding element in $\DC$ is usually denoted by $\dot \ga$.
\end{conv}

Certain structural aspects of definable subsets of $\gC^{\csn, \bullet}$ that are not covered in \cite{Yin:special:trans, Yin:int:acvf, Yin:int:expan:acvf} will be needed below. However, they are of a technical nature and, in order to maintain the central narrative of the paper, we shall defer all such discussions to the last section.

Recall from \cite[Definition~4.6]{Yin:int:expan:acvf} the notion of $\VF$-dimension and its operator $\dim_{\VF}$ for arbitrary definable subsets and from \cite[Definition~3.1]{Yin:int:expan:acvf} the notion of $\RV$-dimension and its operator $\dim_{\RV}$ for definable subsets in the $\RV$-sort. Recall from \cite[Definition~3.14, Lemma~4.9]{Yin:int:expan:acvf} the various variants $\jcb_{\VF}$, $\jcb_{\K}$, etc., of the Jacobian. There is no need to know the particulars of these definitions here and we shall content ourselves with the vague remark that they are sufficiently well-behaved as expected.

\begin{defn}[$\mVF$-categories]\label{defn:f:VF:cat}
An object of the category $\mVF[k]$ is a definable pair $(A, \omega)$, where
\begin{itemize}
  \item $\pr_{\VF}(A) \sub \VF^k$,
  \item $A$ is of $\VF$-dimension $\leq k$ and $\RV$-fiber dimension $0$,
  \item $\omega \coloneqq (\omega_{\K}, \omega_{\Gamma}) : A \fun \K^{\times} \times \Gamma$ is a function, which is understood as a \emph{volume form} on $A$.
\end{itemize}

Let $\bm A = (A, \omega)$ and $\bm B = (B, \sigma)$ be two objects in $\mVF[k]$. A \emph{$\mVF[k]$-morphism} between $\bm A$ and $\bm B$ is a definable \emph{essential bijection} $F : A \fun B$, that is, a bijection that is defined away from certain definable subsets of $A$, $B$ of $\VF$-dimension $< k$, such that, for every point $x \in \dom(F)$,
\begin{gather*}
\omega_{\Gamma}(x) = \sigma_{\Gamma}(F(x)) + \vv(\jcb_{\VF} F(x)),\\
\omega_{\K}(x) = \sigma_{\K}(F(x)) \cdot (\tbk \circ \rv)(\jcb_{\VF}F(x)).
\end{gather*}
\end{defn}

\begin{defn}[$\mRV$-categories]\label{defn:f:RV:cat}
An object of the category $\mRV[k]$ is a definable triple $(U, f, \omega)$, where
\begin{itemize}
  \item $U$ is a subset in the $\RV$-sort,
  \item $f : U \fun (\RV^{\times})^k$ is a function such that $\dim_{\RV}(U_t) = 0$ for all $t \in (\RV^{\times})^k$,
  \item $\omega \coloneqq (\omega_{\K}, \omega_{\Gamma}) : U \fun \K^{\times} \times \Gamma$ is a function, which is understood as a \emph{volume form} on $U$.
\end{itemize}

Let $\bm U = (U, f, \omega)$ and $\bm V = (V, g, \sigma)$ be two objects in $\mRV[k]$. Let $F : U \fun V$ be a definable bijection and
\[
F^{\rightleftharpoons} = \{(t, s) \in f(U) \times g(V) : \ex{ u \in U} (f(u) =  t \wedge (g \circ F)( u) = s) \},
\]
which is understood as the correspondence between $f(U)$ and $g(V)$ induced by $F$. Then $F$ is a \emph{$\mRV[k]$-morphism} between $\bm U$ and $\bm V$ if
\begin{itemize}
  \item for all $(f(u), (g \circ F)(u)) \in F^{\rightleftharpoons}$,
  \[
  \omega_{\Gamma}(u) = \sigma_{\Gamma}(F (u)) + \jcb_{\Gamma} F^{\rightleftharpoons}(f(u), (g \circ F)(u)),
  \]
  \item for all $(f(u), (g \circ F)(u)) \in F^{\rightleftharpoons}$ away from a subset of $F^{\rightleftharpoons}$ of $\RV$-dimension $< k$,
  \[
  \omega_{\K}(u) = \sigma_{\K}(F (u)) \jcb_{\K} F^{\rightleftharpoons}(f(u), (g \circ F)(u)).
  \]
\end{itemize}
\end{defn}

\begin{defn}[$\mRES$-categories]\label{defn:RES:cat}
The category $\mRES[k]$ is the full subcategory of $\mRV[k]$ such that $(U, f, \omega) \in \mRES[k]$ if and only if all coordinates of $U$, $f(U)$ are in $\K$ and $\omega_{\Gamma} = 0$.
\end{defn}

\begin{defn}[$\mG$-categories]\label{def:Ga:cat}
An object of the category $\mG[k]$ is a definable triple $(I, f, \omega)$, where
\begin{itemize}
  \item $I$ is a subset in the $\Gamma$-sort,
  \item $f : I \fun \Gamma^k$ is a function,
  \item $\omega : I \fun \Gamma$ is a function, which is understood as a \emph{volume form} on $I$.
\end{itemize}

Let $\bm I = (I, f, \omega)$ and $\bm J = (J, g, \sigma)$ be two objects in $\mG[k]$. A \emph{$\mG[k]$-morphism} between $\bm I$ and $\bm J$ is a definable bijection $F : I \fun J$ such that, for all $\gamma \in I$,
\[
\omega(\gamma) = \sigma(\gamma) + \jcb_{\Gamma} F^{\rightleftharpoons}(f(\gamma), (g \circ F)(\gamma)).
\]
\end{defn}

Set $\mVF[*] = \coprod_k \mVF[k]$, similarly for $\mRV[*]$, $\mRES[*]$, and $\mG[*]$.

\begin{nota}\label{nota:RV:ele}
Let $(\RV^{\circ \circ})^{\times} = \rv(\MM \mi \{0\})$. We introduce the following shorthand for some elements of the Grothendieck semigroups and their groupifications (and closely related constructions):
\begin{gather*}
[1_{\mu}]_1 = [(\{1\}, \id, \id)] \in \gsk \mRES[1], \quad [0_{\mu}]_1 = [(\{0\}, \id, \id)] \in \gsk \mG[1],\\
[\mathbf{H}_{\mu}]_1 = [((0, \infty), \id, 0)] \in \gsk \mG[1],\\
\mathbf{j}_{\mu} = [((\RV^{\circ \circ})^{\times}, \id, (1, 0))] -  [1_{\mu}]_1 \in \ggk \mRV[1],\\
\mathbf{A}_{\mu} = [(\K^{\times}, \id, (1, 0))] +  [1_{\mu}]_1 \in \ggk \mRES[1].
\end{gather*}
\end{nota}

By \cite[Corollary~3.22]{Yin:int:expan:acvf}, the graded ring $\ggk \KRCV[*]$ may be canonically identified with certain tensor product of the graded rings $\ggk \mRES[*]$, $\ggk \mG[*]$. The following theorem sums up the main results of \cite{Yin:int:expan:acvf} and provides the integration formalism that will be used throughout the rest of the paper.

\begin{thm}[{\cite[Theorems~5.17, 5.21]{Yin:int:expan:acvf}}]\label{main:prop:k:vol:dag}
For each $k \geq 0$ there is a canonical isomorphism of Grothendieck semigroups
\[
\int_{+} \gsk  \mVF[k] \fun \gsk  \KRCV[k] /  \misp,
\]
where $\misp$ is a semigroup congruence relation, such that
\[
\int_{+} [\bm A] = [\bm U]/  \misp \quad \text{if and only if} \quad  [\bm A] = [\bb L \bm U],
\]
where $\bb L : \ob \KRCV[k] \fun \ob \mVF[k]$ is a canonical ``lifting'' map (not defined for the morphisms). Putting these together, we obtain a canonical isomorphism of Grothendieck semirings
\[
\int_{+} \gsk  \mVF[*] \fun \gsk  \KRCV[*] /  \misp.
\]
After groupification, we obtain a canonical graded ring isomorphism
\[
\int \ggk  \mVF[*] \fun \ggk  \KRCV[*] /  \mu \mathbf{I} \coloneqq \bigoplus_{k} \ggk \KRCV[k] / \mu \mathbf{I}_k,
\]
where the ideal $\mu \mathbf{I}$ is generated by the element $\mathbf{j}_{\mu}$, and two graded ring homomorphisms
\begin{gather*}
\Xint{\textup{e}}^g : \ggk \mVF[*] \fun  \ggk \KRCES[*] / (\mathbf{A}_{\mu}),\\
\Xint{\textup{e}}^b: \ggk \mVF[*] \fun \ggk \KRCES[*] / ([1_{\mu}]_1).
\end{gather*}
\end{thm}

\begin{nota}
Henceforth let us abbreviate $\gsk \KRCV[k] / \misp$, $\gsk  \KRCV[*] /  \misp$ as $\gsk \RV_k$, $\gsk \RV$, respectively. Their groupifications are abbreviated accordingly.
\end{nota}

Let $A \sub \VF^n$ be a definable subset. Recall from the discussion after \cite[Theorem~5.17]{Yin:int:expan:acvf} that the set of all definable functions $A \fun \gsk \RV$ is denoted by $\fn(A, \gsk \RV)$, which is a $\gsk \RV$-semimodule. This notation makes sense even when $A$ contains imaginary elements.

\begin{nota}\label{nota:vol}
For any $ \gamma \in \Gamma(S)$ we write $\bm o^{\omega}_{ \gamma}$, $\bm c^{\omega}_{ \gamma}$, and $\bm a^{\omega}_{ \gamma}$ for the elements $\int [(\go(0,  \gamma), \omega)]$, $\int [(\gc(0,  \gamma), \omega)]$, and $\int [(\vv^{-1}( \gamma), \omega)]$ in $\ggk \RV$ and their canonical images in closely related constructions. If $\omega = (1, 0)$ then it shall be omitted from the notation. Note that $\bm a^{\omega}_{\gamma} = \bm a^{(\omega_{\K}, \omega_{\Gamma} + \gamma)}_{0}$, and similar equations hold for $\bm o^{\omega}_{\gamma}$, $\bm c^{\omega}_{\gamma}$. Also,
\[
\bm o_{0} = [1_{\mu}]_1 = \bm a_0 \otimes [\mathbf{H}_{\mu}]_1 \quad \text{and} \quad \bm c_{0} = \bm o_{0} + \bm a_0 = \bm a_0 \otimes ([\mathbf{H}_{\mu}]_1 + [0_{\mu}]_1).
\]
More generally, by \cite[Proposition~3.21]{Yin:int:expan:acvf}, if $\omega_{\Gamma} = 0$ then there are a definable finite partition $D_i$ of $[0, \infty) \sub \Gamma$ with $D_0 = \{0\}$ and twistbacks $\omega_i : \K^{\times} \fun \K^{\times}$ of $\omega_{\K}$ (see \cite[Definition~2.25]{Yin:int:expan:acvf}) such that
\[
\bm o^{\omega}_{0} = \sum_{i > 0} \bm a_0^{\omega_i} \otimes [(D_i, \id, 0)]_1 \quad \text{and} \quad \bm c^{\omega}_{0} = \bm o^{\omega}_{0} + \bm a_0^{\omega_0}.
\]

Let $\KRC$ be the zeroth graded piece of $\ggk \RV[\sqrt{\bm c_0 \bm o_{0}}^{-1}]$, where $\sqrt{\bm c_0 \bm o_{0}}$ is a (formal) square root of $\bm c_0 \bm o_{0}$. The elements $\bm o_{0}\sqrt{\bm c_0 \bm o_{0}}^{-1}$ and $\bm c_{0}\sqrt{\bm c_0 \bm o_{0}}^{-1}$ in $\KRC$ are frequently denoted by $\bm q^{-1}$ and $\bm q$, respectively. Now, the element $\bm 1$ in $\KRC$ may seem obscure since, unlike in discretely valued fields, it does not stand for the volume of any subset. However, this phenomenon is natural if one works with a divisible value group; see \cite[\S1, Remark~3.23]{hru:kazh:val:ring:2006} for an interpretation.

All elements of the forms $\bm o_{\gamma}$, $\bm c_{\gamma}$ in $\KRC$ are units, which is a direct consequence of the change of variables formula and no further localization at these elements is needed; see Corollary~\ref{dual:vol:cons}.

Let $\tau$ be the volume form on $\VF^{\times}$ given by
\[
a \efun (\tbk \circ \rv(a^{-1}), \vv(a^{-1})) = (\rcsn(\rv(a)) / \rv(a), -\vv(a))
\]
and $\tilde \tau = (\tau_{\K}, 0)$. Note that $\tau$ is nothing but a motivic version of the classical volume form $dx/x$ and $\tilde \tau$ its reduction to the residue field. The forms on $\RV^{\times}$ induced by $\tau$, $\tilde \tau$ are also denoted by $\tau$, $\tilde \tau$. We shall simply write $\bm a_0^{\tau} = \bm a_0^{\tilde \tau}$ as $\bm l$ in $\KRC$.
\end{nota}

There is a natural semiring homomorphism
\[
\gsk \RV \fun \ggk \RV \to^{x \efun x \sqrt{\bm c_0 \bm o_{0}}^{-n}} \ggk \RV[\sqrt{\bm c_0 \bm o_{0}}^{-1}] \fun  \KRC,
\]
where $x$ is a homogeneous element of degree $n$. This is of course not guaranteed to be injective. But a significant part of the structure of $\ggk \RV$ does survive in $\KRC$. We define the $\KRC$-module of definable functions $A \fun \KRC$ as
\[
\fn(A, \KRC) = \fn(A, \gsk \RV) \otimes_{\gsk \RV} \KRC.
\]
The induced homomorphism of $\KRC$-modules is understood as integration whose target ring is no longer graded:
\[
\int_A : \fn(A, \KRC) \fun \KRC.
\]

Functions in $\fn(A, \KRC)$ still admit definable representatives of the form $f : A \fun \mdl P(\RV^m)$. (Actually there are various ways to achieve this, all of which are quite tedious. Anyhow, we choose one of them and fix it for the rest of this paper.) For any $\bm f \in \fn(A, \KRC)$, if $f$ is a definable representative of $\bm f$ then we write $[f] = \bm f$. In this case, the equality $[f]( a) = [f]( b)$ means $f( a) =_{ a,  b} f( b)$ in $\KRC_{ a,  b}$ (recall \cite[Notation~3.25]{Yin:int:expan:acvf}). Similarly, if $f'$ is another representative of $\bm f$ then the equality $[f]( a) = [f']( a)$ means $f( a) =_{ a} f'( a)$ in $\KRC_{ a}$. Obviously each element in $\KRC$ may be treated as a constant function in $\fn(A, \KRC)$. However, in general, we cannot treat every constant function $\bm f$ in $\fn(A, \KRC)$ as an element in $\KRC$. Of course, we may do so if there is a definable representative of $f( a)$ for some representative $f$ of $\bm f$ and some $ a \in A$. In particular, if $A$ contains a definable point then the submodule of constant functions of $\fn(A, \KRC)$ may be identified with $\KRC$. For many interesting definable subsets, for example, definable groups, this is indeed the case.

For any definable bijection $\phi : A \fun B \sub \VF^n$ there is the corresponding Jacobian transformation:
\[
\phi^{\jcb} : \fn(A, \KRC) \fun  \fn(B, \KRC).
\]
The following two crucial properties justify conceptually our calling the map $\int_A$ an integration (see \cite[Theorems~5.22, 5.23]{Yin:int:expan:acvf}). For any $\bm f \in \fn(A, \KRC)$,
\begin{itemize}
 \item Fubini property. For any nonempty subsets $E_1, E_2 \sub [n]$,
\[
\int_{ a \in \pr_{E_1}(A)} \int_{A_a} \bm f = \int_{ a \in \pr_{E_2}(A)} \int_{A_a} \bm f.
\]
\item Change of variables. $\int_{A} \bm f = \int_{B} \phi^{\jcb} (\bm f)$.
\end{itemize}

\begin{defn}\label{defn:bd:032014}
A subset $B \sub \VF^n  \times \RV^m$ is \emph{bounded} if both $\vv(\pr_{\VF}(B))$ and $\vrv(\pr_{\RV}(B))$ are bounded from below.

A function $f : \VF^n \fun \VF^m$ is \emph{bounded} if, for every bounded $A \sub \VF^n$, $f(A)$ is also bounded. A function $\bm f \in \fn(A, \KRC)$ has \emph{bounded support} or is \emph{bounded} if the support $\supp(\bm f)$ of $\bm f$ is bounded. The submodule of $\fn(A, \KRC)$ of bounded functions is denoted by $\fn^{b}(A, \KRC)$.
\end{defn}

\begin{conv}
For any volume form $\omega : A \fun \K^{\times} \times \Gamma$ there is a natural $\KRC$-module automorphism of $\fn(A, \KRC)$, which for simplicity is also denoted by $\omega$. If $\omega$ is induced by a scalar $a \in \VF^{\times}$ then we shall simply write it as $\underline{a}$. Note that $\omega$ is a Jacobian transformation if it is the Jacobian of a definable bijection $A \fun A$. The composition
\[
\fn(A, \KRC) \to^{\omega} \fn(A, \KRC) \to^{\int_A} \KRC
\]
is understood as integration of the functions in $\fn(A, \KRC)$ with respect to $\omega$. We shall almost never mention which volume form on $A$ is being used since the results, with few exceptions, do not depend on the choice, as long as it is fixed and satisfies some obvious conditions (for instance, translation-invariance). Accordingly, the integral of $\bm f$ with respect to $\omega$, which should be written as $\int_{(A, \omega)} \bm f$, shall just be written as $\int_{A} \bm f$, unless it is necessary to specify $\omega$. In fact, for simplicity, most of the results below will be stated with respect to the default volume form $(1, 0)$.

The volume of $A$ is written as $\vol(A) = \int_{A} \bm 1 = \int \bm 1_{A}$, where $\bm 1_{A} \in \fn(\VF^n, \KRC)$ is the characteristic function of $A$. For example, $\vol(\MM) = \bm q^{-1}$ and $\vol(\OO) = \bm q$. Note that the volume of $\UU = \OO \mi \MM$ is $\frac{\bm q^2 - \bm 1}{\bm q}$, which is different from its \emph{$\tau$-volume} $\vol_{\tau}(\UU) = \int_{(\UU, \tau)} \bm 1 = \bm l$ (recall the last sentence of Notation~\ref{nota:vol}).

With pointwise multiplication, $\fn(A, \KRC)$ becomes a commutative ring. For $\bm f, \bm g \in \fn(A, \KRC)$ this is simply denoted by $\bm f \bm g$. However, we adopt the convention that if $A$ is a definable group then the default multiplication of $\fn(A, \KRC)$ is given by convolution: for any \emph{translation-invariant} volume form $\omega$ on $A$ and any $f$, $g \in \fn(A, \KRC)$, the convolution $f *_{\omega} g \in \fn(A, \KRC)$ is given by
\[
(f *_{\omega} g)(a) = \int_{(x \in A, \omega)} f(x) g(ax^{-1}).
\]
Clearly if the group operation of $A$ is a bounded function then $\fn^b(A, \KRC)$ is a subring of $\fn(A, \KRC)$.

For various constructions it is often more convenient to modify $\KRC$ further through some standard algebraic operations, such as localization and taking quotient, and then form the module of definable functions as above. It is cumbersome to introduce a new notation whenever this happens. We shall write all of them simply as $\KRC$ and tacitly assume that some operations have been performed so that the construction in question makes sense.
\end{conv}

\begin{exam}\label{exam:sl2}
Let $G = \msl_2$, which also stands for $\msl_2(\VF)$ if no further qualification is indicated (similarly for subsets of $G$). Let $G = KAN$ be the Iwasawa decomposition, where $K = G(\OO)$, $A$ is the subgroup of diagonal matrices, and $N$ is the subgroup of unipotent upper triangular matrices. The opposite of $N$, that is, the subgroup of unipotent lower triangular matrices, is denoted by $\overline N$. The Borel subgroup $B = AN$ of $G$ contains exactly upper triangular matrices and its opposite $\overline B = A \overline N$ contains exactly lower triangular matrices. The Iwahori subgroup $I$ of $K$ is the pre-image of $B(\K)$ under the residue map $\res : G(\OO) \fun G(\K)$.

We shall also write $A(\OO)$, $B(\K)$, etc., as $A_{\OO}$, $B_{\K}$, etc. The group $A_{\csn(\Gamma)}$ is naturally isomorphic to $\Gamma$ and will be written simply as $A_{\Gamma}$. The quotient $A / A_{\OO}$ is naturally isomorphic to $A_{\OO} \times A_{\Gamma}$, which is also denoted by $\chi_*(A)$ and is understood as the group of cocharacters of $A$. The extended Weyl group $\mathbf{N}_{G}(A) / A_{\OO}$ is denoted by $\tilde W$, where $\mathbf{N}_{G}(A)$ is the normalizer of $A$ in $G$, and the Weyl group $\mathbf{N}_{K}(A) / A_{\OO}$ by $W$, where $\mathbf{N}_{K}(A) = \mathbf{N}_{G}(A) \cap K$. The standard representatives of $W$ are chosen to be $\bigl(\begin{smallmatrix}
1 & 0 \\ 0 & 1
\end{smallmatrix}\bigr)$ and $w \coloneqq \bigl(\begin{smallmatrix}
0 & 1 \\ -1 & 0
\end{smallmatrix}\bigr)$. The subgroup generated by these representatives is denoted by $W^*$. We have $\tilde W = W \ltimes \chi_*(A)$, where the action of $W$ is given by that of the sign on $\Gamma$.

Clearly $\dim_{\VF}(G) = 3$. Let $\wh G = \pr_{< 4}(G) \cap \VF^3$, similarly for $K$, $A$, $N$, etc. For integration over $G$, we may identify $G$ with $\wh G$, since $\dim_{\VF}(G \mi \pr_{< 4}^{-1}(\wh G)) = 2$. Hence all volume forms on $G$ may be thought of as
\[
d \begin{pmatrix}
         x & y\\
         z & w
       \end{pmatrix} = \omega(x, y, z) dx \wedge dy \wedge dz,
\]
where $\omega : \wh G \fun \K^{\times} \times \Gamma$ is a definable function. For example, the unimodular Haar volume form on $G$ is given by $(x, y, z) \efun \tau(x)$ if $x \neq 0$, which we shall also denote by $\tau$. We have $\vol_{\tau}(I) = \bm{l}$. To compute $\vol_{\tau}(K)$, we first write $\wh K = \wh K' \cup \wh K''$, where $\pr_1(\wh K') = \UU$ and $\pr_1(\wh K'') = \MM$. We have $\vol_{\tau}(\wh K') = \bm l \bm q^2$. Since
\[
\wh K'' = \bigcup_{\gamma > 0} \vv^{-1}(\gamma) \times \Big( \bigcup_{a \in \UU} \{a\} \times (a^{-1} + \gc(0, \gamma)) \Big),
\]
we see that
\[
\vol_{\tau}(\wh K'') = \vol_{(\tau_{\K}, 0)} ( \MM \times \UU \times \OO ) = \bm o^{\tilde \tau}_0 (\bm q^2 - \bm 1)
\]
and hence
\[
\vol_{\tau}(K) = \bm l \bm q^2 + \bm o^{\tilde \tau}_0 \bm q^2 - \bm o^{\tilde \tau}_0 = \bm c^{\tilde \tau}_0 \bm q^2 - \bm o^{\tilde \tau}_0.
\]

Let $C_c(I \setminus G /I)$ be the Iwahori-Hecke algebra of the $I$-bi-invariant bounded functions in $\fn^{b}(G, \KRC)$ with the Haar form, where multiplication is given by the normalized convolution:
\[
(\bm f * \bm g)(a) = \bm{l}^{-1} \int_{x \in G} \bm f(x) \bm g(ax^{-1}).
\]
This makes $C_c(I \setminus G /I)$ a unital commutative algebra with the multiplicative identity $\bm 1_I$. The structure of $C_c(I \setminus G /I)$ has been investigated in \cite{hru:kazh:val:ring:2006}, where the ring $\KRC$ is further simplified to be a field in a canonical way. It seems that more can be done following the treatment of \cite{iwo:alg:HPK}.
\end{exam}

\section{Integrable functions}\label{section:integrable}
Recall that $\Omega \coloneqq \VF / \MM$, a subset of imaginary elements, is viewed as a motivic analogue of the subgroup of roots of unity of the unit circle and the quotient map $\theta : \VF \fun \Omega$ a motivic analogue of a generic additive character.

In the classical situation, summation of a character over any compact subgroup vanishes. To accommodate this phenomenon when we integrate with an additive character, a cancellation rule must and shall be introduced in the construction of the ring in which integration with an additive character takes values. Also, in order to avoid undesirable collapses of volumes of definable subsets, we need to work with bounded subsets of $\Omega$.

For simplicity, we shall write $\fn(-, \KRC)$ as $\fn(-)$ when there is no danger of confusion. From here on the technicalities discussed in the last section will be needed.

\begin{lem}\label{gam:fin:O}
Let $O$ be a $\Gamma$-algebraic subset of $\Omega$. Then, for every $\K$-coset $\Delta \sub \Omega$, $O \cap \Delta$ is finite.
\end{lem}
\begin{proof}
This follows easily from \cite[Corollary~3.4]{Yin:int:expan:acvf}.
\end{proof}

\begin{defn}
Let $V$ be a definable abelian group and $W$ a definable subgroup of $V$. A definable function $\bm f \in \fn(V)$ is \emph{$W$-invariant} if, for any $a \in V$ and any $b \in W$, $\bm f(a) = \bm f(a b)$. The subgroup of $\fn(V)$ of $W$-invariant functions is denoted by $\fn(V)^W$.
\end{defn}

\begin{lem}\label{quot:id}
The group $\fn(\VF)^{\MM}$ (resp.\ $\fn(\VF)^{\OO}$, $\fn(\VF^{\times})^{1+\MM}$) is isomorphic to the group $\fn(\Omega)$ (resp.\ $\fn(\VF/ \OO)$, $\fn(\RV)$).
\end{lem}
\begin{proof}
The arguments being the same, we will just prove this for the additive case of $\MM$. It suffices to show that, for each element $\bm f \in \fn(\VF)^{\MM}$ with a representative $f : \VF \fun \mdl P(\RV^m)$, there is a definable function $f_{\downarrow} : \Omega \fun \mdl P(\RV^m)$ such that, for very $\omega \in \Omega$, there is an $a \in \omega$ with $f(a) = f_{\downarrow}(\omega)$. To that end, let $\pi$ be a $\Gamma$-partition of $f$. Each $\pi^{-1}( \gamma)$ may be naturally extended to a function
\[
f_{ \gamma} : \VF \fun \mdl P(\RV^m \times \{1, \infty\}).
\]
By \cite[Lemma~3.18]{Yin:special:trans} (see \cite[Remark~3.19]{Yin:special:trans} for the multiplicative case), for each $f_{ \gamma}$, there are finitely many $\MM$-cosets away from which $f_{ \gamma}$ is a partial function on $\Omega$. Let $\Delta$ be the collection of these exceptional $\MM$-cosets. So $f \rest (\VF \mi \theta^{-1}(\Delta))$ is a partial function on $\Omega$. On the other hand, by Lemmas~\ref{gam:alg} and \ref{di:cen}, $\Delta$ has definable centers. The lemma follows.
\end{proof}

Convolution in $\fn(\Omega)$ is normalized as in Example~\ref{exam:sl2} and consequently $\fn(\Omega)$ becomes a unital commutative ring. By Lemma~\ref{quot:id} we may and shall identify $\fn(\VF)^{\MM}$ with $\fn(\Omega)$. Similarly $\fn(\OO)^{\MM}$, $\fn(\VF)^{\OO}$ may be identified with $\fn(\K)$, $\fn(\VF/\OO)$, respectively. We are interested in the subrings $\fn^b(\Omega)$, $\fn(\K)$ and the subgroup $\fn^b(\VF/\OO)$ of bounded functions (see Definition~\ref{defn:bd:032014}). Of course $\fn(\K)$ is also a subring of $\fn^b(\Omega)$. On the other hand, $\fn^b(\VF / \OO)$ may be naturally identified with a subgroup of $\fn^b(\Omega)$. Since $\fn^b(\VF / \OO)$ does not contain units and, for $\bm f \in \fn^b(\Omega)$ and $\bm g \in \fn^b(\VF / \OO)$,
\[
\bm f * \bm g \in \fn^b(\VF / \OO),
\]
we see that $\fn^b(\VF / \OO)$ is actually an ideal of $\fn^b(\Omega)$.

Let $\fn^b(\Omega)_{\fin}$ be the subring of $\fn^b(\Omega)$ consisting of those definable functions whose support is covered by a (bounded) $\Gamma$-algebraic subset of $\K$-cosets $\Delta \sub \Omega$ (here $\Delta$ is identified with $\dot \Delta \in \DC$).

\begin{lem}\label{fun:omega:decom}
For every $\bm f \in \fn^b(\Omega)$ there are an $\bm f' \in \fn^b(\Omega)_{\fin}$ and an $\bm f'' \in \fn^b(\VF / \OO)$ such that $\bm f = \bm f' + \bm f''$.
\end{lem}
\begin{proof}
This is immediate by (the proof of) Lemma~\ref{quot:id}.
\end{proof}

Next we consider the ring $\fn(\K \times \Omega)$, where convolution is again normalized as in Example~\ref{exam:sl2}, which is compatible with convolution in $\fn(\Omega)$ via the Fubini property for $\KRC$-valued integrals. It follows from Lemma~\ref{gam:alg} that the definable functions in $\fn(\K \times \Omega)$ with the property that the projection of its support into $\Omega$ is contained in a bounded $\Gamma$-algebraic subset form a subring. This subring may be identified with the ring of definable functions $\Omega \fun \fn(\K)$ with support contained in a bounded $\Gamma$-algebraic subset. Although it is not the group ring $\fn(\K)[\Omega]$ in the usual sense, we abuse the notation slightly to denote it as such. Also note that if we write a typical element $\bm p$ in $\fn(\K)[\Omega]$ as a formal sum $\sum_{\omega \in O} \omega \bm f_{\omega}$, where $O$ is a $\Gamma$-algebraic subset of $\Omega$, then any representative of $\bm f_{\omega}$ is $\omega$-definable (but may not be definable). We call $\bm f_{\omega}$ the coefficient of $\omega$. If each $\bm f_{\omega}$ is a constant function then $\bm p$ is a \emph{$\K$-constant} function. Let $\fn^1(\K)[\Omega]$ be the subgroup of $\fn(\K)[\Omega]$ of $\K$-constant functions. Since multiplication is given by convolution in $\fn(\K \times \Omega)$, it is easy to see that $\fn^1(\K)[\Omega]$ is an ideal of $\fn(\K)[\Omega]$.

Since $\K$ is a subgroup of $\Omega$, there are two natural actions of $\K$ on $\fn(\K \times \Omega)$: one on the $\K$-coordinate and the other on the $\Omega$-coordinate. We identify them as follows. For any $\bm f, \bm g \in \fn(\K \times \Omega)$, if there are $t, t' \in \K$ such that, for every $(s, \omega) \in \K \times \Omega$,
\[
\bm g(s, \omega) = \bm f(t + s, t' + \omega)
\]
then we write $A_{(t,t')}(\bm f) = \bm g$. The \emph{anti-diagonal ideal} $\mdl{A}$ of $\fn(\K \times \Omega)$ is generated by elements of the form $\bm f - A_{(t, -t)}(\bm f)$. For simplicity, the ideal $\mdl{A} \cap \fn(\K)[\Omega]$ of $\fn(\K)[\Omega]$ is also denoted by $\mdl A$.

\begin{thm}[{\cite[Theorem~11.3]{hrushovski:kazhdan:integration:vf}}]\label{cancel:rule}
There is a canonical ring isomorphism
\[
\phi :  \fn^b(\Omega) / \fn^b(\VF / \OO) \fun \fn(\K)[\Omega] / (\fn^1(\K)[\Omega] + \mdl A).
\]
\end{thm}
\begin{proof}
For $\omega \in \Omega$ and $\bm f \in \fn^b(\Omega)$ we denote the translation of $\bm f$ by $\omega$ as $T_{\omega}(\bm f)$; that is, the $\omega$-definable function $T_{\omega}(\bm f)$ is given by $T_{\omega}(\bm f)(t) = \bm f(t + \omega)$. Let
\[
\rho : \fn(\K)[\Omega] \fun \fn^b(\Omega)_{\fin}
\]
be the ring homomorphism given by
\[
\sum_{\omega \in O} \omega \bm f_{\omega} \efun \sum_{\omega \in O} T_{-\omega} (\bm f_{\omega}),
\]
where $O \sub \Omega$ is a bounded $\Gamma$-algebraic subset and, by Lemma~\ref{gam:fin:O}, the formal sum on the right-hand side may be naturally interpreted as an element in $\fn^b(\Omega)_{\fin}$. We claim that $\rho$ is surjective. To see this, let $\bm f \in \fn^b(\Omega)_{\fin}$ and $D$ be a bounded $\Gamma$-algebraic subset of $\K$-cosets $\Delta \sub \Omega$ that covers the support of $\bm f$. By Lemma~\ref{di:cen}, there is a definable subset $A \sub \VF$ such that $A \cap \bigcup \Delta$ is a singleton $\{a_{\Delta}\}$ for every $\dot \Delta \in D$. Let $\omega_{\Delta} = a_{\Delta} + \MM$. So $\{\omega_{\Delta} : \dot \Delta \in D \} \sub \Omega$ is a bounded $\Gamma$-algebraic subset. For each $\omega_{\Delta}$ let $\bm f_{\omega_{\Delta}} = T_{\omega_{\Delta}}(\bm f) \rest \K$. Clearly
\[
\rho \Big(\sum_{\dot \Delta \in D} \omega_{\Delta} \bm f_{\omega_{\Delta}} \Big) = \bm f.
\]
On the other hand, by Lemma~\ref{fun:omega:decom}, there is a canonical surjective homomorphism
\[
\sigma: \fn^b(\Omega)_{\fin} \fun \fn^b(\Omega) / \fn^b(\VF / \OO).
\]
It is not hard to see that the kernel of the surjective homomorphism $\sigma \circ \rho$ is precisely the ideal $\fn^1(\K)[\Omega] + \mdl A$.
\end{proof}

The quotient ring
\[
\fn(\K)[\Omega] / (\fn^1(\K)[\Omega] + \mdl A)
\]
may be regarded as a motivic analogue of the group ring $\R[S^1]$, where $S^1$ is the unit circle on the complex plane, and is henceforth denoted by $\KCC$. The elements of $\Omega$ in $\KCC$ are treated as symbols rather than $\MM$-cosets. The conceptual reason that we do not directly take the quotient ring
\[
\fn^b(\Omega) / \fn^b(\VF / \OO)
\]
as the analogue is that it is not completely free of $\VF$-data, where elements of $\Omega$ and $\VF/\OO$ are treated as subsets of $\VF$. However, technically, it is more effective to work with the latter, which is what we shall do below, since its construction corresponds directly to the cancellation rule mentioned above.

Multiplication in $\KCC$ comes from the normalized convolution in $\fn^b(\Omega)$. For each $\bm x \in \KRC$ let $e_{\bm x} \in \fn^b(\Omega)$ be the function given by $\omega \efun \bm x$ if $\omega = 0$ and $\omega \efun 0$ otherwise. Clearly if $\bm x \neq \bm x'$ then $e_{\bm x} - e_{\bm x'} \notin \fn^b(\VF / \OO)$. So $\KRC$ may be treated as a subring of $\KCC$ via the canonical embedding
\[
\bm x \efun e_{\bm x} + \fn^b(\VF / \OO).
\]
Similarly, if for each $\omega \in \Omega$ we let $e_{\omega} \in \fn^b(\Omega)$ be the function with support $\set{\omega}$ and $e_{\omega}(\omega) =  1$, then the mapping given by
\[
\omega \efun e_{\omega} + \fn^b(\VF / \OO)
\]
is an injective group homomorphism from $\Omega$ into the multiplicative group of $\KCC$. So $\Omega$ may be considered as a multiplicative subgroup of $\KCC$.

\begin{nota}
In writing we shall not distinguish $\KRC$ and $\Omega$ from their canonical images in $\fn^b(\Omega)$ or $\KCC$. Also, to emphasize the transition from an additively written group to a multiplicatively written group, we shall think of the embedding of $\Omega$ into $\fn^b(\Omega)$ or $\KCC$ as an exponential map and denote it by $\exp$.
\end{nota}

For any definable subset $A \sub \VF^n$ we can again form the $\KCC$-module $\fn(A, \KCC)$ of definable functions. Clearly $\fn(A, \KCC)$ contains $\fn(A)$ as a subgroup (or a sub-$\KRC$-module). Similarly $\fn^b(A, \KCC)$ contains $\fn^b(A)$ as a subgroup.

\begin{defn}\label{def:integrable}
A function $\bm f \in \fn(A, \KCC)$ is \emph{integrable} if there is a representative $\tilde{\bm f} : A \fun \fn^b(\Omega)$ of $\bm f$ that is uniformly bounded, that is, there is a $\gamma \in \Gamma$ such that $\supp(\tilde{\bm f}( a)) \sub \go(0, \gamma)$ for every $ a \in A$. In this case, we say that $\tilde{\bm f}$ is an \emph{integrable representative} of $f$. It is easy to see that $\bm f$ is integrable if and only if there is a definable $\gamma  \in \Gamma$ such that, for every representative $\tilde{\bm f} : A \fun \fn^b(\Omega)$ of $\bm f$ and every $ a \in A$, $\tilde{\bm f}( a)$ is $\OO$-invariant outside $\go(0, \gamma)$.

We say that $\bm f$ is \emph{almost integrable} if $\bm f \rest (A \cap \go(0, \gamma))$ is integrable for every $\gamma \in \Gamma$. By compactness, $\bm f$ is almost integrable if and only if there is a representative $\tilde{\bm f}$ of $\bm f$ and a definable function $h : \Gamma \fun \Gamma$ such that $\tilde{\bm f} \rest (A \cap \go(0, \gamma))$ is uniformly bounded by $h(\gamma)$.

If for every $ a \in A$ there is an $ a$-definable $\gamma \in \Gamma$ such that $\bm f \rest (A \cap \go( a, \gamma))$ is integrable then $\bm f$ is \emph{locally integrable}.
\end{defn}

\begin{nota}\label{nota:int}
The subsets of $\fn(A, \KCC)$ of integrable functions, almost integrable functions, and locally integrable functions are denoted by $\ifn(A)$, $\ifn^{a}(A)$, and $\ifn^{l}(A)$, respectively. If $A = \VF^n$ then we simply write $\ifn_n$, etc. Obviously these are subgroups of $\fn(A, \KCC)$ and are closed under pointwise multiplication. We have
\[
\ifn(A) \subsetneq \ifn^a(A) \subsetneq \ifn^l(A).
\]
\end{nota}

Let $\bm g : A \fun \fn(\Omega)$ be a definable function. For each $\omega \in \Omega$, the function $\leb_{\omega}(\bm g) : A \fun \KRC_{\omega}$ is given by $ a \efun \bm g( a)(\omega)$. Let $\leb(\bm g) \in \fn(\Omega)$ be the function given by $\omega \efun \int_{A} \leb_{\omega}(\bm g)$.

\begin{prop}\label{integration:C}
There is a canonical $\KCC$-module homomorphism
\[
\int_A : \ifn(A) \fun \KCC.
\]
\end{prop}
\begin{proof}
Let $\tilde{\bm f}, \tilde{\bm f}' : A \fun \fn^b(\Omega)$ be two integrable representatives of an integrable function $\bm f \in \ifn(A)$. For any $\K$-coset $\Delta$ and any $\omega, \rho \in \Delta$, since
\[
\leb_{\omega} (\tilde{\bm f})( a) - \leb_{\omega} (\tilde{\bm f}')( a) \quad \text{and} \quad \leb_{\rho} (\tilde{\bm f})( a) - \leb_{\rho} (\tilde{\bm f}')( a)
\]
can both be represented by the same $( a, \dot \Delta)$-definable subset, where $\dot \Delta$ is the imaginary element corresponding to $\bigcup \Delta$, we have
\begin{align*}
\int_{A} \leb_{\omega} (\tilde{\bm f}) - \int_{A} \leb_{\omega} (\tilde{\bm f}') &= \int_{A} (\leb_{\omega} (\tilde{\bm f}) - \leb_{\omega} (\tilde{\bm f}')) \\
&= \int_{A} (\leb_{\rho} (\tilde{\bm f}) - \leb_{\rho} (\tilde{\bm f}')) \\
&= \int_{A} \leb_{\rho} (\tilde{\bm f}) - \int_{A} \leb_{\rho} (\tilde{\bm f}'),
\end{align*}
which is an element in $\KRC_{\dot \Delta}$. So $\leb (\tilde{\bm f}) - \leb (\tilde{\bm f}') \in \fn^b(\VF / \OO)$. So if we set
\[
\int_A \bm f = \leb (\tilde{\bm f}) +  \fn^b(\VF / \OO)
\]
then $\int_A$ is a well-defined $\KCC$-module homomorphism.
\end{proof}

\begin{defn}
Let $E \sub [n]$. We say that a definable function $\bm f \in \fn(A, \KCC)$ is \emph{iteratively integrable on $E$} if
\begin{itemize}
  \item $\bm f \rest A_a \in \ifn(A_a)$ for every $ a \in \pr_E (A)$,
  \item $\pr_E(\bm f) \in \ifn(\pr_{E}(A))$, where $\pr_E(\bm f)$ is the function given by $ a \efun \int_{A_a} \bm f$.
\end{itemize}
In other words, $\bm f$ is iteratively integrable on $E$ if and only if the iterated integral
\[
\int_{ x \in \pr_{E} (A)} \int_{ y \in A_x} \bm f( x,  y)
\]
is defined. It clearly follows from Definition~\ref{def:integrable} and compactness that this iterated integral is defined if and only if there is a representative $\tilde{\bm f} : A \fun \fn^b(\Omega)$ of $\bm f$ and a $\gamma \in \Gamma$ such that, for every $ a \in \pr_E (A)$,
\begin{itemize}
 \item $\tilde{\bm f} \rest A_a$ is uniformly bounded,
 \item $\leb (\tilde{\bm f} \rest A_a)$ is $\OO$-invariant outside $\go(0, \gamma)$.
\end{itemize}
Let $\ifn_E(A)$ denote the group of functions that are iteratively integrable on $E$.
\end{defn}

Obviously $\ifn(A) \sub \ifn_E(A)$ for any $E$. By the definition of $\int_A$ in Proposition~\ref{integration:C} and the Fubini property of $\KRC$-valued integrals, we also have:

\begin{prop}\label{fubini}
For all $E_1, E_2 \sub [n]$ and every $\bm f \in \ifn_{E_1}(A) \cap \ifn_{E_2}(A)$,
\[
\int_{\pr_{E_1} (A)} \pr_{E_1}(\bm f) = \int_{\pr_{E_2} (A)} \pr_{E_2} (\bm f).
\]
\end{prop}

Therefore, the Fubini property also holds for integrable functions. It is also easy to see that the change of variables formula holds for integrable functions.

\begin{lem}\label{inv:para}
Let $\phi : A \fun \VF^{m} \times \Gamma^l$ be a definable function. For any $\bm f, \bm g \in \ifn(A)$, if $\int_{\phi^{-1}( a,  \gamma)} \bm f = \int_{\phi^{-1}( a,  \gamma)} \bm g$ for all $( a,  \gamma)$ then  $\int_{A} \bm f = \int_{A} \bm g$.
\end{lem}
\begin{proof}
This is immediate by Proposition~\ref{integration:C} and compactness.
\end{proof}

\begin{defn}
Suppose that $A$ is a subgroup of $\VF^n$. A function $f \in \fn(A, \KCC)$ is a \emph{definable additive character} of $A$ if it is a group homomorphism $A \fun \Omega$.

For any $ b \in \VF^n$ let $\chi_{ b}$ be the $ b$-definable map given by $ a \efun \theta( a \cdot  b)$, where $ a \cdot  b$ is the ordinary dot product. Clearly $\exp_{ b} = \exp \circ \chi_{ b}$ is an additive character of $\VF^n$. Note that $\exp_{ b}$ is almost integrable but not integrable.
\end{defn}

The notation $\gb^{\iota}$ is explained in Definition~\ref{defn:dual}.

\begin{lem}\label{simp:char:int}
Suppose that $t \in \Omega$ is definable. For any $ a,  b \in \VF^n$ and any polydisc $\gb$ around $ a$ with $\dim_{\VF}(\gb) = n$,
\[
\int_{ x \in \gb} \exp(\chi_{ b}( x) + t) = \begin{dcases}
          \vol(\gb) \exp(\chi_{ b}( a) + t), &\text{if }  b \in \gb^{\iota} -  a;\\
          \bm 0, &\text{otherwise}.
\end{dcases}
\]
\end{lem}
\begin{proof}
To simplify the argument, we shall just show the case that $\gb$ is an open polydisc $\go( a, \gamma)$. The proof for the other cases are almost identical. For the first equality, since $ b \in \gc(0, - \gamma)$, clearly $\chi_{ b}( a') = \chi_{ b}( a)$ for every $ a' \in \gb$. So the equality is clear by the definition of $\int_{\gb}$.

For the second equality, we first consider the case $n = 1$ and hence $ a$, $ b$ are simply written as $a$, $b$. Observe that, since $\vv(b) < -\gamma$ and $b \gb = \go(ba, \vv(b) + \gamma)$, the image $\chi_b(\gb)$ is a union of $\K$-cosets. We have $(\omega - t) / b \sub \gb$ for any $\omega - t \in \chi_b(\gb)$. Note that, for any $a' \in \gb$, $\chi_b(a') + t = \omega$ if and only if $a' \in (\omega - t) / b$. Let $\bm h : \gb \fun \fn^b(\Omega)$ be an integrable representative of $(\chi_b + t) \rest \gb$. Then it is enough to show that $\leb(\bm h) \in \fn^b(\VF/\OO)$. In fact, $\leb(\bm h)$ is a constant function on its support $\chi_b(\gb) + t$. To see this, observe that, for any $d \in \VF$ with $\theta(d) \in \chi_b(\gb) + t$,
\[
\leb(\bm h(\theta(d))) = \int_{\gb} \leb_{\theta(d)} (\bm h) = \vol ((\theta(d) - t)/b) \in \KRC_d.
\]
By Lemma~\ref{di:cen}, $t$ contains a definable point. So, for any other $d' \in \VF$, the obvious $(d, d')$-definable bijection between $(\theta(d) - t)/b$ and $(\theta(d') - t)/b$ implies that
\[
\vol((\theta(d) - t)/b) =_{d, d'} \vol((\theta(d') - t)/b).
\]

We now consider the case $n > 1$. Without loss of generality, we may assume $\vv(b_1) < -\gamma$. Let $ b_1 = (b_2, \ldots, b_n)$. For each $ a_1 \in \pr_{> 1} (\gb)$, by the case $n = 1$ above, we have
\[
\pr_{> 1} (\exp (\chi_{ b} + t))( a_1) = \int_{x_1 \in \pr_1(\gb)} \exp(\chi_{b_1}(x_1) + \theta( a_1 \cdot  b_1) + t) = \bm 0.
\]
By the Fubini property,
\[
\int_{ x \in \gb} \exp(\chi_{ b}( x) + t) = \int_{ y \in \pr_{> 1}(\gb)} \pr_{> 1} (\exp(\chi_{ b} + t))( y) = \bm 0,
\]
as desired.
\end{proof}

Recall the various special volumes from Notation~\ref{nota:vol}.

\begin{lem}[Averaging formula]\label{averaging}
Let $\bm f \in \ifn(A)$ and $p : A \fun \Gamma$ be an $\go$-partition such that $\go( a, p( a)) \sub A$ for every $ a \in A$. Then
\[
\int_{A} \bm f = \int_{ x \in A} \bm o_{p( x)}^{-n} \int_{ y \in \go( x, p( x))} \bm f( y).
\]
\end{lem}
\begin{proof}
Let $W = \bigcup_{ a \in A} (\set{ a} \times \go( a, p( a)))$ and $\bm f^*$ be the definable function on $W$ given by $( a,  b) \efun \bm o_{p( a)}^{-n} \bm f( b)$. Note that this is well-defined since, by Corollary~\ref{dual:vol:cons}, $\bm o_{p( a)}$ is invertible in $\KCC_{ a}$. It is easy to see that any integrable representative of $\bm f$ gives rise to an integrable representative of $\bm f^*$ and hence $\bm f^*$ is integrable. So, by the Fubini property, the righthand side of the above equation is well-defined.

Now, by Lemma~\ref{inv:para}, it is enough to show that, for any $\gamma \in \Gamma$,
\[
\int_{p^{-1}(\gamma)} \bm f = \bm o_{\gamma}^{-n} \int_{ x \in p^{-1}(\gamma)} \int_{ y \in \go( x, \gamma)} \bm f( y).
\]
Let $\bm g$ be the function on $p^{-1}(\gamma) \times \go(0, \gamma)$ such that $\bm g( a,  b) = \bm f^*( a,  a +  b)$. Clearly $\bm g$ is integrable and, by change of variables,
\[
\bm o_{\gamma}^{-n} \int_{ x \in p^{-1}(\gamma)} \int_{ y \in \go( x, \gamma)} \bm f( y) = \int_{ x \in p^{-1}(\gamma)} \int_{ y \in \go(0, \gamma)} \bm g( x,  y).
\]
Notice that, for every $ b \in \go(0, \gamma)$,
\[
\int_{ x \in p^{-1}(\gamma)} \bm g( x,  b) = \int_{ x \in p^{-1}(\gamma)} \bm f^*( x,  x +  b) = \bm o_{\gamma}^{-n} \int_{ x \in p^{-1}(\gamma)} \bm f( x +  b) = \bm o_{\gamma}^{-n} \int_{p^{-1}(\gamma)} \bm f,
\]
where the last equality is by change of variables again. Finally, by the Fubini property, we have
\[
\int_{ x \in p^{-1}(\gamma)} \int_{ y \in \go(0, \gamma)} \bm g( x,  y) = \int_{ y \in \go(0, \gamma)} \int_{ x \in p^{-1}(\gamma)} \bm g ( x,  y) = \bm o_{\gamma}^n \bm o_{\gamma}^{-n} \int_{p^{-1}(\gamma)} \bm f = \int_{p^{-1}(\gamma)} \bm f,
\]
as desired.
\end{proof}

\begin{prop}\label{convol:C}
Suppose that $A$ is a definable subgroup of $\VF^n$. The convolution product of two integrable functions on $A$ always exists and is also an integrable function. Moreover, the convolution map
\[
* : \ifn(A)^2 \fun \ifn(A)
\]
is $\KCC$-bilinear, associative, and commutative.
\end{prop}
\begin{proof}
Let $\bm f, \bm g \in \ifn(A)$ and $\tilde{\bm f}$, $\tilde{\bm g}$ be two integrable representatives of $\bm f$, $\bm g$ that are uniformly bounded by $\alpha, \beta \in \Gamma$, respectively. Let $\tilde{\bm h}_{ a}$ be the function on $A$ given by $ b \efun \tilde{\bm f}( b) \tilde{\bm g}( a - b)$. Clearly the support of every $\tilde{\bm h}_{ a}( b)$ is bounded by $\min \set{\alpha, \beta}$ and hence the support of $\leb (\tilde{\bm h}_{ a})$ is bounded by $\min \set{\alpha, \beta}$. So the function given by $ a \efun \leb (\tilde{\bm h}_{ a})$ is an integrable representative of $\bm f * \bm g$.

From the definition of convolution, $\KCC$-bilinearity is clear. For associativity, given a third $\bm h \in \ifn(A)$, we have
\begin{align*}
((\bm f * \bm g) * \bm h)( a) &= \int_{ y \in A} \Big( \int_{ x \in A} \bm f( x) \bm g( y -  x) \Big) \bm h( a -  y) &\\
& = \int_{ x \in A} \bm f( x)  \int_{ y \in A} \bm  g( y -  x) \bm  h( a -  y)  & &\text{by the Fubini property}\\
& = \int_{ x \in A} \bm f( x) \int_{ z \in A}  \bm g( z) \bm h( a -  z -  x)  & &\text{by change of variables}\\
& = (\bm f * (\bm g * \bm h))( a).
\end{align*}
Commutativity may be proved in a similar way, since the standard proof also only makes use of the Fubini property and change of variables.
\end{proof}

Set $\ifn^b(A) = \ifn(A) \cap \fn^b(A, \KCC)$. Suppose that $A$ is a definable subgroup of $\VF^n$. For any two functions $\bm f, \bm g \in \ifn^b(A)$ and any $ a,  c \in A$ with $\vv( a)$ sufficiently low, either $\bm f( c) = \bm 0$ or $\bm g( a -  c) = \bm 0$, and hence $\int_{ y \in A} \bm f( y) \bm g( a - y) = \bm 0$. This means:

\begin{cor}\label{convol:closed:bounded}
The group $\ifn^b(A)$ is closed under convolution and hence is a commutative ring.
\end{cor}

\section{Integration with an additive character}\label{section:with:add:char}

In this section we shall define the Fourier transform for definable functions on $\VF^n$. We mention two conceptual reasons why this definition should work (and it does). The first is self-duality of local fields, that is, the group of additive characters of any local field $L$ may be identified with the additive group of $L$ itself. This fact makes integrating over the group of definable characters possible in our first-order setting. The second is that we have included in $\KCC$ a tautological image $\Omega$ of $\VF$ via the generic additive character $\exp$ of $\VF$. It is conceivable that the whole construction may be adapted for any definable abelian group $G$ as long as the ``dual group'' $G^*$ of $G$, whatever that may mean, can be handled in a first-order fashion (for example, if $G$ is an abelian variety) and the ring $\KCC$ is so enlarged that it contains the images of both $G$ and $G^*$ under definable characters.

\begin{defn}
The \emph{Fourier transform} of a function $\bm f \in \fn(\VF^n, \KCC)$ is the function $\wh{\bm f} \in \fn(\VF^n, \KCC)$ such that
\[
\wh{\bm f}( b) = \int_{ x \in \VF^n} \bm f( x) \exp_{ b}( x),
\]
provided that the integral is defined for every $ b \in \VF^n$. Sometimes $\wh{\bm f}$ is also written as $\mdl F(\bm f)$. This is more suggestive if we understand the Fourier transform as a (partial) linear operator.
\end{defn}

Recall from Notation~\ref{nota:int} that $\ifn^b_n$ is short for $\ifn^b(\VF^n)$, and so on.

\begin{lem}\label{fourier:exists}
If $\bm f \in \ifn^b_n$ then the integral $\int_{ x \in \VF^n} \bm f( x) \exp_{ b}( x)$ is defined for every $ b \in \VF^n$.
\end{lem}
\begin{proof}
Let $\bm f \in \ifn^b_n$ and $\tilde{\bm f} : \VF^n \fun \fn^b(\Omega)$ be an integrable representative of $\bm f$. Let $\go(0, \alpha)$ be a definable open polydisc that contains both $\supp(\leb (\tilde{\bm f}))$ and $\supp(\bm f)$. For any $ b \in \VF^n$, clearly there is a $\beta \in \Gamma$, which depends on $\alpha$ and $\vv( b)$, such that $\supp( \tilde{\bm f}( a) \exp_{ b}( a) ) \sub \go(0, \beta)$
for every $ a \in \go(0, \alpha)$. So the integral in question is defined.
\end{proof}

\begin{rem}\label{bound:supp:inb}
The above lemma shows that the Fourier transform of any $\bm f \in \ifn^b_n$ exists. On the other hand, since we need to consider $\vv( b)$, there is no guarantee that $\wh{\bm f}$ is integrable. Of course, if $\wh{\bm f}$ has bounded support then it is integrable. In other words, $\wh{\bm f}$ is always an almost integrable function. The map $\mdl F : \ifn^b_n \fun \ifn^a_n$ is clearly a $\KCC$-module homomorphism.
\end{rem}

\begin{prop}[Convolution formula]\label{convol:formula:fun}
Let $\bm f, \bm g \in \ifn^b_n$. Then
\[
\mdl F(\bm f * \bm g) = \mdl F(\bm f) \mdl F(\bm g).
\]
\end{prop}
\begin{proof}
Note that, by Corollary~\ref{convol:closed:bounded}, $\bm f * \bm g \in \ifn^b_n$ and hence the transform $\mdl F(\bm f * \bm g)$ exists. For any definable $\gamma \in \Gamma$ that is sufficiently low, by the Fubini property and change of variables, we have
\begin{align*}
\mdl F(\bm f * \bm g)( b) &=  \int_{ x \in \go(0,\gamma)}  \int_{ y \in \go(0,\gamma)} \bm f( y) \bm g( x - y) \exp_{ b}( x)\\
 &=  \int_{ y \in \go(0,\gamma)}  \int_{ z \in \go(0, \gamma)} \bm f( y) \bm g( z)  \exp_{ b}( z +  y) \\
 &= \Big( \int_{ y \in \go(0,\gamma)} \bm f( y) \exp_{ b}( y) \Big) \Big( \int_{ z \in \go(0, \gamma)} \bm g( z)  \exp_{ b}( z) \Big) \\
 & = \wh{\bm f}( b) \wh{\bm g}( b),
\end{align*}
as desired.
\end{proof}

Let $A \sub \VF^n$ be a definable subset. A function $\bm f \in \fn(A, \KCC)$ is \emph{locally constant at $ a \in A$} if there is a $\gamma \in \Gamma$ such that $\bm f \rest (\go( a, \gamma) \cap A)$ is constant. It is \emph{locally constant} if it is locally constant at every $ a \in A$. Obviously if $A$ is discrete then every function in $\fn(A, \KCC)$ is locally constant. If $\bm f$ is locally constant at $ a \in A$ then, by compactness, \cite[Lemma~2.21]{Yin:int:expan:acvf}, and \omin-minimality, there is an $ a$-definable $h( a) \in \Gamma$ such that $\bm f \rest (\go( a, h( a)) \cap A)$ is constant.

\begin{defn}
If $\bm f$ is locally constant then we associate with it a definable function $\iota_{\bm f} : A \fun \Gamma$ as follows. By compactness, there is a definable function $h : A \fun \Gamma$ such that $\bm f \rest (\go( a, h( a)) \cap A)$ is constant for every $ a \in A$. Let
\[
B_{ a} = \{ b \in A :  a \in \go( b, h( b)) \text{ and } \vv( a) \leq h( b) \leq h( a) \}.
\]
By \omin-minimality, there are finitely many $ a$-definable points in $\Gamma$ that determine $h(B_{ a})$. Let $\iota_{\bm f}( a)$ be the lowest of these points. It is clear that $\iota_{\bm f}( a) = \iota_{\bm f}( a')$ for any $ a' \in \go( a, \iota_{\bm f}( a)) \cap A$. That is, $\iota_{\bm f}$ is an $\go$-partition of $A$.

If $\bm f$ has bounded support then there is a definable open polydisc $\go(0, \alpha)$ that contains $\supp(\bm f)$. In this case, we take $\iota_{\bm f}( a) = \vv( a)$ for every $ a \in A \mi \go(0, \alpha)$. If $\bm f$ is the characteristic function of a definable subset $V \sub A$ then we write $\iota_V( a)$ if $\bm f$ is locally constant at $ a$ and $\iota_V : A \fun \Gamma$ if $\bm f$ is locally constant.

Note that, although the construction of $\iota_{\bm f}$ depends on the choice of $h$, the discussion below does not depend on this choice.
\end{defn}

\begin{thm}[Fourier inversion formula]\label{fourier:inverse}
For any $\bm f \in \ifn^b_n$, if $\wh{\bm f} \in \ifn^b_n$ then, for almost all $a \in \VF^n$, that is, for all $a \in \VF^n$ away from a definable subset of $\VF^n$ of $\VF$-dimension $< n$,
\[
\dhat{\bm f}(- a) = \bm f( a).
\]
\end{thm}
\begin{proof}
Applying Lemma~\ref{fun:loc:cons} to any definable representative of $\bm f$, we see that there is a definable open subset $A \sub \VF^n$ with $\dim_{\VF}(\VF^n \mi A) < n$ such that $\bm f \rest A$ is locally constant. Fix an $ a \in A$ and let $\alpha = \iota_{\bm f \rest A}( a)$. For any definable open polydisc $\go(0, \gamma)$ that contains $ a$ and the supports of $\bm f$ and $\wh{\bm f}$, by the Fubini property, we have
\begin{align*}
\dhat{\bm f}(- a) &= \int_{ x \in \VF^n} \Big(  \int_{ y \in \VF^n} \bm f( y) \exp_{ x}( y) \Big) \exp_{-  a}( x) \\
&= \int_{ y \in \go(0,\gamma)}  \int_{ x \in \go(0,\gamma)} \bm f( y) \exp_{ y -  a}( x).
\end{align*}
If $\gamma$ is so low that $\go(0, \gamma)^{\iota} \sub \go(0, \alpha)$ then, by Lemma~\ref{simp:char:int} and change of variables,
\[
\int_{ y \in  a + (\go(0, \gamma) \mi \go(0, \gamma)^{\iota})} \int_{ x \in \go(0,\gamma)}  \bm f( y) \exp_{ y -  a}( x) = \bm 0,
\]
and hence
\[
\dhat{\bm f}(- a) = \int_{ y \in  a + \go(0, \gamma)^{\iota}}  \int_{ x \in \go(0,\gamma)} \bm f( y) \exp_{ y -  a}( x) = \bm c_{\gamma}^n \bm o_{\gamma}^n \bm f( a) = \bm f( a),
\]
where the last equality is by Corollary~\ref{dual:vol:cons}.
\end{proof}

Now, as in the classical integration theory, there is a submodule of $\ifn^b_n$ whose image under $\mdl F$ is contained in $\ifn^b_n$ and hence the double Fourier transform of any function in the submodule exists.

\begin{defn}
The submodule of $\ifn^b_n$ consisting of Schwartz-Bruhat functions, that is, locally constant integrable functions with bounded support, is called the \emph{Schwartz space on $\VF^n$} and is denoted by $\mathscr S_n$.
\end{defn}

\begin{lem}\label{iota:bound:exist}
For any $\bm f \in \mathscr S_n$, $\iota_{\bm f}(\VF^n)$ is bounded from above.
\end{lem}
\begin{proof}
Let $\go(0, \gamma)$ be a definable open polydisc containing $\supp(\bm f)$ such that if $ a \notin \go(0, \gamma)$ then $\iota_{\bm f}( a) = \vv( a) \leq \gamma$. Since $\iota_{\bm f} \rest \go(0, \gamma)$ is an $\go$-partition of $\go(0, \gamma)$, by Lemma~\ref{vol:par:bounded}, $\iota_{\bm f}(\go(0, \gamma))$ is bounded from above.
\end{proof}

From now on we shall need Notation~\ref{nota:vol:m}.

\begin{lem}\label{con:vol}
Let $a \in \VF^{\times}$ be definable. For any $\bm f \in \mathscr S_n$,
\[
\int_{(A, \underline{a})} \bm f = \bm m^n_{\vv(a)} \int_{A} \bm f.
\]
\end{lem}
\begin{proof}
Let $\beta \in \Gamma$ be a definable element such that $\beta \pm \vv(a) > \iota_{\bm f}(\VF^n)$. By the averaging formula and the change of variables formula,
\[
\int_{(A, \underline{a})} \bm f = \int_{ x \in A}  \bm o_{\beta}^{-n} \int_{(\go(0, \beta), \underline{a})} \bm f( x) = \int_{ x \in A}  \bm o_{\beta}^{-n} \int_{\go(0, \beta + \vv(a))} \bm f( x) = \bm m^n_{\vv(a)} \int_{A} \bm f,
\]
as required.
\end{proof}

\begin{prop}\label{sb:tran:sb}
The $\KCC$-module $\mathscr S_n$ is closed under Fourier transform.
\end{prop}
\begin{proof}
Fix a nonzero $\bm f \in \mathscr S_n$. By Lemma~\ref{iota:bound:exist} and \omin-minimality, there is a definable upper bound $\beta \in \Gamma$ of $\iota_{\bm f}(\VF^n)$. By Lemma~\ref{averaging}, for any $ b \in \VF^n$ we have
\begin{align*}
\wh{\bm f}( b)  &= \int_{ x \in \VF^n} \Big( \bm o_{\iota_{\bm f} ( x)}^{-n} \int_{ y \in \go( x, \iota_{\bm f} ( x))} \bm f( y) \exp_{ b}( y) \Big)  \\
&= \int_{ x \in \VF^n} \Big( \bm f( x) \bm o_{\iota_{\bm f} ( x)}^{-n} \int_{ y \in \go( x, \iota_{\bm f} ( x))} \exp_{ b}( y) \Big).
\end{align*}
If $ b \notin \gc(0, -\beta)$ then, by Lemma~\ref{simp:char:int}, for any $ a \in \VF^n$, $\int_{\go( a, \iota_{\bm f} ( a))} \exp_{ b} = 0$ and hence $\wh{\bm f}( b) = 0$. So $\wh{\bm f}$ has bounded support. By Remark~\ref{bound:supp:inb}, $\wh{\bm f}$ is integrable.

Choose a definable $\alpha \in \Gamma$ such that $\gc(0, \alpha)$ contains $\supp(\bm f)$. Let $ b \in \VF^n$ and $ b' \in \go( b, -\alpha)$. If $ a \notin \supp(\bm f)$ then obviously
\[
\bm f( a) \exp_{ b}( a) = \bm f( a) \exp_{ b'}( a) = \bm 0.
\]
If $ a \in \supp(\bm f)$ then, for any $ a' \in \go( a, \iota_{\bm f} ( a))$, since $\go( a, \iota_{\bm f} ( a)) \sub \gc(0, \alpha)$, we have
\[
 b \cdot  a' -  b' \cdot  a' = ( b -  b') \cdot  a' \in \MM
\]
and hence $\int_{\go( a, \iota_{\bm f} ( a))} \exp_{ b} = \int_{\go( a, \iota_{\bm f} ( a))} \exp_{ b'}$. So, for every $ a \in \VF^n$,
\[
\bm f( a) \bm o_{\iota_{\bm f} ( a)}^{-n} \int_{\go( a, \iota_{\bm f} ( a))} \exp_{ b} = \bm f( a) \bm o_{\iota_{\bm f} ( a)}^{-n} \int_{\go( a, \iota_{\bm f} ( a))} \exp_{ b'}.
\]
So $\wh{\bm f}$ is also locally constant.
\end{proof}

\begin{cor}
The restriction of $\mdl F$ is a $\KCC$-module automorphism of $\mathscr S_n$.
\end{cor}
\begin{proof}
It is a $\KCC$-module homomorphism by Proposition~\ref{sb:tran:sb}. By the Fourier inversion formula, it must be a bijection.
\end{proof}

For any function $\bm f \in \fn(\VF^n, \KCC)$, $\Check{\bm f}$ is the function given by $\bm f(-  a) = \Check{\bm f}( a)$.

\begin{cor}\label{inverse:convol}
For any $\bm f, \bm g \in \mathscr S_n$,
\[
\mdl F(\bm f \bm g) = \mdl F(\bm f) * \mdl F(\bm g).
\]
\end{cor}
\begin{proof}
By the convolution formula, the Fourier inversion formula, and the fact that $\mathscr S_n$ is closed under convolution,
\begin{align*}
(\mdl F(\bm f) * \mdl F(\bm g))( a) &= \mdl F(\mdl F(\mdl F(\bm f) * \mdl F(\bm g)))(-  a)\\
  &= \mdl F(\mdl F(\mdl F(\bm f)) \mdl F(\mdl F(\bm g)))(-  a)\\
   &=\mdl F(\Check{\bm f} \Check{\bm g})(-  a)\\
    &= \mdl F(\bm f \bm g)( a),
\end{align*}
as desired.
\end{proof}

\begin{thm}[Plancherel formula]\label{plancherel}
Suppose that $\bm f, \wh{\bm f}, \bm g \in \ifn^b_n$. Then
\[
\int_{\VF^n} \bm f \bm g = \int_{\VF^n} \wh{\bm f} \wh{\Check{\bm g}}.
\]
\end{thm}
\begin{proof}
First note that, by Remark~\ref{bound:supp:inb}, $\wh{\bm f} \wh{\Check{\bm g}} \in \ifn^b_n$ and hence the righthand side of the equality is well-defined. For any definable $\gamma \in \Gamma$ that is sufficiently low,
\begin{align*}
 \int_{ x \in \VF^n} \wh{\bm f}( x) \wh{\Check{\bm g}}( x) &=   \int_{ x \in \go(0, \gamma)} \wh{\bm f}( x) \int_{ z \in \go(0, \gamma)} \bm g( z) \exp_{ x}(-  z) \\
  &=   \int_{ z \in \go(0, \gamma)} \bm g( z) \int_{ x \in \go(0, \gamma)} \wh{\bm f}( x) \exp_{-  z}( x) \\
  &= \int_{ z \in \go(0, \gamma)} \bm g( z) \bm f( z),
\end{align*}
where the third equality is by the Fourier inversion formula.
\end{proof}

We note that, for any $\bm f, \bm g \in \ifn^b_n$, a straightforward computation using only the Fubini property shows that the following classical version of the Plancherel formula holds:
\[
\int_{\VF^n} \mdl F(\bm f) \bm g = \int_{\VF^n} \bm f \mdl F ( \bm g).
\]

\section{Definable distributions}\label{section:dist}

Our main references for the classical theory of distributions are \cite{gelfand:shilov:1964, hormander:83}.

By a definable function $\VF^n \times \Gamma \fun \KCC$ we mean a function $\bm f \in \fn(\VF^{n+1}, \KCC)$ such that, for every $( a, b_1)$, $( a, b_2) \in \VF^{n+1}$ with $\vv(b_1) = \vv(b_2)$, $\bm f( a, b_1) = \bm f( a, b_2)$. We write $\fn(\VF^n \times \Gamma, \KCC)$, $\ifn_{n, \Gamma}$, etc., for the corresponding $\KCC$-modules of definable functions. If $\bm f \in \fn(\VF^n \times \Gamma, \KCC)$ then the function $\bm f_{\gamma} : \VF^n \fun \KCC$ is given by $ a \efun \bm f( a, \gamma)$ and the function $\bm f_{ a} : \Gamma \fun \KCC$ is given by $\gamma \efun \bm f( a, \gamma)$. By Lemma~\ref{quot:id}, each $\bm f_{ a}$ may be treated as an $ a$-definable function on $\RV$.

\begin{defn}[{\cite[Definition~11.5]{hrushovski:kazhdan:integration:vf}}]
A \emph{predistribution} on $\VF^n$ is a definable function $\gD \in \fn(\VF^n \times \Gamma, \KCC)$ such that
\begin{itemize}
  \item $\gD_{\gamma}$ is an almost integrable function for every $\gamma \in \Gamma$,
  \item for every $( a, \gamma), ( a', \gamma) \in \VF^n \times \Gamma$, if $\vv( a -  a') > \gamma$, that is, if $\go( a, \gamma) = \go( a', \gamma)$, then $\gD( a, \gamma) = \gD( a', \gamma)$,
 \item $\gD$ is \emph{coherent}, that is, for every $( a, \gamma), ( a, \gamma') \in \VF^n \times \Gamma$, if $\gamma' \geq \gamma$ then
 \[
  \gD( a, \gamma) = \int_{ x \in \go( a, \gamma)} \bm o_{\gamma'}^{-n} \gD( x, \gamma').
 \]
\end{itemize}
If $\bm f \in \ifn^a_n$ then, as in the classical theory of distributions, the almost integrable function $\gD_{\bm f}$ on $\VF^n \times \Gamma$ given by
\[
\gD_{\bm f}( a, \gamma) = \int_{\go( a, \gamma)} \bm f
\]
is a predistribution on $\VF^n$, which shall be called a \emph{regular} predistribution. That $\gD_{\bm f}$ is coherent is clear by the averaging formula (Lemma~\ref{averaging}).
\end{defn}

\begin{lem}
Let $\gD_1$, $\gD_2$ be two predistributions. Suppose that for every $ a \in \VF^n$ there is a $\gamma_{ a} \in \Gamma$ such that $\gD_1( b, \gamma) = \gD_2( b, \gamma)$ for every $\gamma \geq \gamma_{ a}$ and every $ b \in \go( a, \gamma_{ a})$. Then $\gD_1 = \gD_2$.
\end{lem}
\begin{proof}
For each $ a \in \VF^n$ let $G_{ a} \sub \Gamma$ be the $ a$-definable subset of all the values that satisfy the given property. Suppose for contradiction that there is a $( d, \alpha)$ such that $\gD_1( d, \alpha) \neq \gD_2( d, \alpha)$. Then every $G_{ a}$ is bounded from below. By coherence of $\gD_1$, $\gD_2$ and \omin-minimality, there is a least element $\beta_{ a}$ in every $G_{ a}$. Let $p : \VF^n \fun \Gamma$ be the definable function given by $ a \efun \beta_{ a}$. Note that if $ b \in \go( a, p( a))$ then $p( b) = p( a)$ and hence $p$ is an $\go$-partition of $\VF^n$. So $p \rest \go( d, \alpha)$ is an $\go$-partition of $\go( d, \alpha)$. By Lemma~\ref{vol:par:bounded}, there is a $( d, \alpha)$-definable $\alpha' \in \Gamma$ such that $p(\go( d, \alpha)) < \alpha'$. Since $\alpha < p( d) < \alpha'$, by coherence of $\gD_1$, $\gD_2$,
\[
\gD_1( d, \alpha) = \int_{ x \in \go( d, \alpha)} \bm o_{\alpha'}^{-n} \gD_1( x, \alpha') = \int_{ x \in \go( d, \alpha)} \bm o_{\alpha'}^{-n} \gD_2( x, \alpha') = \gD_2( d, \alpha),
\]
which is a contradiction.
\end{proof}

\begin{lem}\label{pre:dis:vol:par}
Let $\gD$ be a predistribution and $\bm f \in \ifn^a_n$ a locally constant function. Suppose that $p: \go( a, \gamma) \fun \Gamma$ is a definable $\go$-partition such that $p( b) \geq \iota_{\bm f}( b) \geq \gamma$ for every $ b \in \go( a, \gamma)$. Then
\[
\int_{ x \in \go( a, \gamma)} \bm f( x) \bm o^{-n}_{p( x)} \gD( x, p( x)) \quad \text{and} \quad \int_{ x \in \go( a, \gamma)} \bm f( x) \bm o^{-n}_{\iota_{\bm f}( x)} \gD( x, \iota_{\bm f}( x))
\]
are defined and they are equal.
\end{lem}
\begin{proof}
Let $\beta \in \Gamma$ be an $( a, \gamma)$-definable bound of $p(\go( a, \gamma)) \cup \iota_{\bm f}(\go( a, \gamma))$, which exists by Lemma~\ref{vol:par:bounded}. For any $ b \in \go( a, \gamma)$, by coherence,
\[
\gD( b, p( b)) = \int_{ y \in \go( b, p( b))} \bm o_{\beta}^{-n} \gD( y, \beta) \quad \text{and} \quad \gD( b, \iota_{\bm f}( b)) = \int_{ y \in \go( b, \iota_{\bm f}( b))} \bm o_{\beta}^{-n} \gD( y, \beta).
\]
Therefore,
\begin{align*}
\bm o_{\beta}^{-n} \int_{ x \in \go( a, \gamma)} \bm f( x) \gD( x, \beta) &= \bm o_{\beta}^{-n} \int_{ x \in \go( a, \gamma)}\bm o^{-n}_{p( x)} \int_{ y \in \go( x, p( x))}  f( y) \gD( y, \beta)\\
&= \int_{ x \in \go( a, \gamma)} \bm f( x) \bm o^{-n}_{p( x)} \int_{ y \in \go( x, p( x))} \bm o_{\beta}^{-n} \gD( y, \beta)\\
&= \int_{ x \in \go( a, \gamma)} \bm f( x) \bm o^{-n}_{p( x)} \gD( x, p( x)),
\end{align*}
where the first equality is by the averaging formula. A completely similar computation shows
\[
\bm o_{\beta}^{-n} \int_{ x \in \go( a, \gamma)} \bm f( x) \gD( x, \beta) = \int_{ x \in \go( a, \gamma)} \bm f( x) \bm o^{-n}_{\iota_{\bm f}( x)} \gD( x, \iota_{\bm f}( x)).
\]
The lemma follows.
\end{proof}

\begin{cor}\label{pre:dis:var:rad}
Let $\gD$ be a predistribution and $p : \go( a, \gamma) \fun \Gamma$ an $\go$-partition. Then
\[
\gD( a, \gamma) = \int_{ x \in \go( a, \gamma)} \bm o_{p( x)}^{-n} \gD( x, p( x)).
\]
Therefore, if $p$ is an $\go$-partition of $\VF^n$ then the function given by $ b \efun \bm o_{p( b)}^{-n} \gD( b, p( b))$ is an almost integrable function.
\end{cor}

By Lemma~\ref{pre:dis:vol:par}, the following is well-defined:

\begin{defn}\label{defn:distr}
Let $\gD$ be a predistribution. A \emph{distribution induced by $\gD$} is a linear functional $\gD' : \mathscr{S}_n \fun \KCC$ given by
\[
\gD'(\bm f) = \int_{ x \in \VF^n} \bm f( x) \bm o^{-n}_{\beta} \gD( x, \beta),
\]
where $\beta$ is any definable upper bound of $\iota_{\bm f}(\VF^n)$, which exists by Lemma~\ref{iota:bound:exist}.

If $\gD$ is regular then $\gD'$ is a \emph{regular} distribution; otherwise $\gD'$ is a \emph{singular} distribution. If $\bm h \in \ifn^a_n$ then, by the averaging formula, we have $\gD_{\bm h} (\bm f) = \int_{ x \in \VF^n} \bm f( x) \bm h( x)$. So the regular distribution $\gD_{\bm h}$ is naturally defined on a larger domain, namely $\ifn^b_n$; that is, $\gD_{\bm h} : \ifn^b_n \fun \KCC$ is a linear functional given by the above integral.
\end{defn}

For each $ a \in \VF^n$, $\gD'$ induces an $ a$-definable function $\gD'_{ a} : \Gamma \fun \KCC$ given by $\gamma \efun \gD'(\bm 1_{\go( a, \gamma)})$. It is clear that $\gD'_{ a}(\gamma) = \gD( a, \gamma)$ for all $( a, \gamma) \in \VF^n \times \Gamma$. Therefore, there is no need to distinguish notationally a predistribution from the distribution induced by it.

\begin{lem}[{\cite[Proposition~11.8]{hrushovski:kazhdan:integration:vf}}]\label{distribution:appox:fun}
Let $\gD$ be a distribution. Then there is a definable open subset $V \sub \VF^n$ and a locally constant function $\bm f \in \ifn^a(V)$ such that $\dim_{\VF}(\VF^n \mi V) < n$ and, for every bounded clopen definable subset $U \sub V$,
\[
\gD(\bm 1_U) = \int_{U} \bm f.
\]
\end{lem}
\begin{proof}
By Lemma~\ref{fun:loc:cons} and \omin-minimality, there is a definable open subset $V \sub \VF^n$ with $\dim_{\VF}(\VF^n \mi V) < n$ and an $\go$-partition $p$ of $V$ such that, for all $ a,  b \in V$, if $ b \in \go( a, p( a))$ then $\gD_{ a} = \gD_{ b}$. For any $ a \in V$ and any $\beta \geq \alpha \geq p( a)$, by coherence, we have
\[
\gD( a, \alpha) = \int_{ x \in \go( a, \alpha)} \bm o^{-n}_{\beta} \gD( x, \beta) = \int_{ x \in \go( a, \alpha)} \bm o^{-n}_{\beta} \gD( a, \beta) = \bm o^n_{\alpha} \bm o^{-n}_{\beta} \gD( a, \beta)
\]
and hence
\[
\bm o^{-n}_{\alpha}\gD( a, \alpha) = \bm o^{-n}_{\beta} \gD( a, \beta).
\]
Let $\bm f$ be the definable function on $V$ given by $ a \efun \bm o^{-n}_{p( a)} \gD( a, p( a))$, which, by Corollary~\ref{pre:dis:var:rad}, is an almost integrable locally constant function such that $\iota_{\bm f} = p$. Let $U \sub V$ be a bounded clopen definable subset and $\beta \in \Gamma$ a definable upper bound of $\iota_U(\VF^n) \cup p(U)$, which exists by Lemma~\ref{vol:par:bounded}. Then
\[
\gD(\bm 1_U) = \int_{ x \in U} \bm o^{-n}_{\beta} \gD( x, \beta) = \int_{ x \in U} \bm f(x),
\]
as required.
\end{proof}

\begin{rem}
An analogue of Bernstein's theorem in~\cite{bernstein:joseph:1972} for nonarchimedean local fields is proved in~\cite[Corollary~11.10]{hrushovski:kazhdan:integration:vf}. It says the following: For any integers $n$, $d$, any nonarchimedean local field $L$ of sufficiently large residue characteristic, and any polynomial $G \in L[X_1, \ldots, X_n]$ of degree $\leq d$, there is a proper variety $V_G \sub L^n$ such that, away from $V_G$, the Fourier transform $\mdl F(\abs{G})$ agrees with a locally constant function. Of course $G$ may be replaced by any definable function $\VF^n \fun \VF$. This is a direct consequence of Lemma~\ref{distribution:appox:fun} and the specialization procedure described in \cite[\S6]{Yin:int:expan:acvf}. A detailed presentation of this will appear elsewhere.
\end{rem}

A simple computation shows that the Fourier transform $\wh{\bm 1}_{\go( a, \gamma)}$ may be written as $\bm o_{\gamma}^{n} \exp_{ a} \bm 1_{\gc(0, - \gamma)}$. By Proposition~\ref{sb:tran:sb}, $\wh{\bm 1}_{\go( a, \gamma)} \in \mathscr S_n$. If $\vv( a) > \gamma$ then we may assume $\iota_{\wh{\bm 1}_{\go( a, \gamma)}}( b) = - \gamma$ for all $ b \in \gc(0, - \gamma)$; if $\vv( a) \leq \gamma$ then we may assume $\iota_{\wh{\bm 1}_{\go( a, \gamma)}}( b) = - \vv( a)$ for all $ b \in \gc(0, - \gamma)$.

\begin{lem}\label{pre:dis:fourier}
For any distribution $\gD$, the function $\wh \gD$ on $\VF^n \times \Gamma$ given by $\wh{\gD} ( a, \gamma) = \gD(\wh{\bm 1}_{\go( a, \gamma)})$ is a predistribution.
\end{lem}
\begin{proof}
Let $\beta \geq - \gamma$ be $\gamma$-definable. Then, for every $ a \in \VF^n$,
\[
\gD(\wh{\bm 1}_{\go( a, \gamma)}) = \bm o_{\gamma}^{n} \int_{ x \in \gc(0, - \gamma)} \exp_{ a}( x) \bm o_{\beta}^{-n} \gD( x, \beta).
\]
Since $\gD_{\beta}$ is an almost integrable function and the function given by $( y,  x) \efun \exp_{ y}( x)$ is also almost integrable, clearly $\wh{\gD}_{\gamma}$ is an almost integrable function as well. Next, if $\go( a, \gamma) = \go( a', \gamma)$ then $\exp_{ a} \rest \gc(0, - \gamma) = \exp_{ a'} \rest \gc(0, - \gamma)$ and hence $\gD(\wh{\bm 1}_{\go( a, \gamma)}) = \gD(\wh{\bm 1}_{\go( a', \gamma)})$. Lastly, let $\gamma' \geq \gamma$. Then
\begin{align*}
\wh{\gD}( a, \gamma) &= \int_{ y \in \gc(0, - \gamma)} \bm o_{\beta}^{-n} \gD( y, \beta) \bm o_{\gamma}^{n} \exp_{ a}( y) \\
&= \int_{ y \in \gc(0, - \gamma')} \bm o_{\beta}^{-n} \gD( y, \beta) \int_{ x \in \go( a, \gamma)} \exp_{ x}( y) & &\text{by Lemma~\ref{simp:char:int}}\\
&= \int_{ x \in \go( a, \gamma)} \bm o_{\gamma'}^{-n} \bm o_{\gamma'}^{n} \int_{ y \in \gc(0, - \gamma')} \exp_{ x}( y) \bm o_{\beta}^{-n} \gD( y, \beta) \\
&= \int_{ x \in \go( a, \gamma)} \bm o_{\gamma'}^{-n} \wh{\gD}( x, \gamma'),
\end{align*}
where the last equality holds since $\beta \geq - \gamma'$.
\end{proof}

\begin{defn}
The \emph{Fourier transform} of the distribution induced by the predistribution $\gD$ is the distribution induced by the predistribution $\wh{\gD}$. Sometimes $\wh{\gD}$ is also written as $\mdl F(\gD)$.
\end{defn}

\begin{thm}\label{dis:fourier}
For every $\bm f \in \mathscr{S}_n$, $\gD(\wh{\bm f}) = \wh{\gD}(\bm f)$.
\end{thm}
\begin{proof}
By Proposition~\ref{sb:tran:sb}, $\wh{\bm f} \in \mathscr S_n$. Let $\beta$ be a definable upper bound of $\iota_{\bm f}(\VF^n)$ such that $\supp(\wh{\bm f}) \sub \gc(0, - \beta)$. Let $\delta \geq - \beta$ be a definable upper bound of $\iota_{\wh{\bm f}}(\VF^n)$. Then, for every definable $\gamma \in \Gamma$ that is sufficiently low,
\begin{align*}
\wh{\gD}(\bm f) &= \int_{ x \in \go(0, \gamma)} \bm f( x) \bm o_{\beta}^{-n} \bm o_{\beta}^{n} \int_{ y \in \gc(0, - \beta)} \exp_{ x}( y) \bm o_{\delta}^{-n} \gD( y, \delta) \\
&= \int_{ y \in \gc(0, - \beta)} \Big( \int_{ x \in \go(0, \gamma)} \bm f( x)  \exp_{ y}( x) \Big) \bm o_{\delta}^{-n} \gD( y, \delta),
\end{align*}
which is just $\gD(\wh{\bm f})$.
\end{proof}

\begin{defn}
Let $\gD$ be a distribution on $\VF^n$. We say that $\gD$ \emph{vanishes at $ a \in \VF^n$} if there is a $\gamma \in \Gamma$ such that $\gD( a', \gamma') = \bm 0$ for every $ a' \in \go( a, \gamma)$ and every $\gamma' \geq \gamma$. By Definition~\ref{defn:distr}, if $\gD$ vanishes at $ a$ then $\gD(\bm f) = 0$ for every $\bm f \in \mathscr{S}_n$ with $\supp(\bm f) \sub \go( a, \gamma)$. If $\gD$ does not vanish at $ a$ then $ a$ is an \emph{essential point} of $\gD$.

The \emph{support} of $\gD$ is the subset $\supp(\gD) \sub \VF^n$ of essential points of $\gD$.
\end{defn}

\begin{rem}\label{opar:dist}
In general, $\supp(\gD)$ is \emph{not} guaranteed to be definable. However, it is not hard to see that, by compactness, if $\supp(\gD)$ is bounded then there is a definable $\gamma \in \Gamma$ such that $\supp(\gD) \sub \gc(0,\gamma)$.

By Lemma~\ref{fun:loc:cons} and \omin-minimality, if $\gD$ vanishes at $ a$ then there is an $ a$-definable $\iota_{\gD}( a) \in \Gamma$ such that $\gD( a', \gamma') = \bm 0$ for every $ a' \in \go( a, \iota_{\gD}( a))$ and every $\gamma' \geq \iota_{\gD}( a)$. Moreover, if $ a' \in \go( a, \iota_{\gD}( a))$ then $\gD$ vanishes at $ a'$ and $\iota_{\gD}( a') = \iota_{\gD}( a)$. Therefore, for any parametrically definable subset $A \sub \VF^n \mi \supp(\gD)$, by compactness, $\iota_{\gD}$ may be constructed as a parametrically definable $\go$-partition of $A$. If, in addition, $A$ is bounded and closed then, by Lemma~\ref{vol:par:bounded}, there is a parametrically definable upper bound of $\iota_{\gD}(A)$. Also, if $\supp(\gD) \sub \gc(0,\gamma)$ then we may set $\iota_{\gD}( a) = \vv( a)$ for every $ a \notin \gc(0, \gamma)$.
\end{rem}

\begin{lem}\label{dis:supp:bound}
For any $\go( a, \gamma) \sub \VF^n \mi \supp(\gD)$, $\gD( a, \gamma) = 0$. In particular, $\gD$ has no essential point if and only if $\gD = \bm 0$.
\end{lem}
\begin{proof}
The first assertion follows from Remark~\ref{opar:dist} and coherence. The second assertion is the special case $\supp(\gD) = \0$.
\end{proof}

\begin{lem}\label{dis:bou:four:reg}
If $\gD$ has bounded support then there is a locally constant function $\bm h \in \ifn^a_n$ such that $\wh{\gD} = \gD_{\bm h}$.
\end{lem}
\begin{proof}
Fix a definable $\gamma \in \Gamma$ such that $\gc(0, \gamma)$ contains $\supp(\gD)$. Let $\bm h$ be the definable function on $\VF^n$ given by
\[
 a \efun \gD(\exp_{ a} \rest \gc(0, \gamma)) = \int_{ x \in \gc(0, \gamma)} \exp_{ a}( x) \bm o^{-n}_{\beta_{ a}} \gD( x, \beta_{ a}),
\]
where $\beta_{ a} = \max \set{\pm \gamma}$ if $ a \in \go(0, - \gamma)$ and $\beta_{ a} = \max \set{\pm \gamma, \pm \vv( a)}$ if $ a \notin \go(0, - \gamma)$. It is clear that $\bm h$ is an almost integrable locally constant function and $\iota_{\bm h}( a) = \beta_{ a}$. Now it is enough to show that, for every $\bm f \in \mathscr{S}_{n}$,
\[
\gD(\wh{\bm f}) = \gD_{\bm h} (\bm f) = \int_{ x \in \VF^n} \bm f( x) \bm h( x).
\]
Observe that, for every definable $\alpha \in \Gamma$ that is sufficiently low and every definable $\beta \geq - \alpha$, $\bm h \rest \go(0, \alpha)$ is given by
\[
 a \efun \int_{ x \in \gc(0, \gamma)} \exp_{ a}( x) \bm o^{-n}_{\beta} \gD( x, \beta)
\]
and hence, by the Fubini property,
\begin{align*}
\gD_{\bm h} (\bm f) &= \int_{ y \in \go(0, \alpha)} \bm f( y) \int_{ x \in \gc(0, \gamma)} \exp_{ y}( x) \bm o^{-n}_{\beta} \gD( x, \beta)\\
&= \int_{ x \in \gc(0, \gamma)} \wh{\bm f}( x) \bm o^{-n}_{\beta} \gD( x, \beta).
\end{align*}
By Lemma~\ref{dis:supp:bound}, for every $ a \notin \gc(0, \gamma)$, $\gD( a, \beta) = 0$. Since $\beta \geq \iota_{\wh{\bm f}}(\VF^n)$, we deduce that $\gD_{\bm h} (\bm f) = \gD(\wh{\bm f})$.
\end{proof}

\begin{defn}
For any $\bm h \in \ifn^a_n$, let $\bm h \gD$ and $\gD \bm h$ denote the linear functional $\mathscr{S}_n \fun \KCC$ given by
\[
\bm f \efun (\gD\bm h)( \bm f) = \int_{ x \in \VF^n} \bm f ( x) \bm o^{-n}_{\beta} \gD( x, \beta) \bm h ( x),
\]
where $\beta$ is any definable upper bound of $\iota_{\bm h}(\VF^n)$. Since $\bm f \bm h \in \ifn^b_n$, this is well-defined by Lemma~\ref{pre:dis:vol:par}. By the averaging formula, this integral may also be written as
\[
\int_{ x \in \VF^n} \bm o^{-2n}_{\beta} \bm f ( x) \gD( x, \beta) \int_{ y \in \go( x, \beta)} \bm h ( y).
\]
Therefore, if $\gD_{\bm h} = \gD_{\bm g}$ then $\bm h \gD = \bm g \gD$. This means that, for any regular distribution $\gD'$, it makes sense to form their product $\gD \gD' = \gD' \gD$, where the distribution is given by $\gD \bm h$ for any $\bm h \in \ifn^a_n$ such that $\gD' = \gD_{\bm h}$.
\end{defn}

\begin{defn}
Let $\gD_1$, $\gD_2$ be predistributions on $\VF^{n_1}$, $\VF^{n_2}$, respectively. Let $\gD_1 \otimes \gD_2$ be the definable function on $\VF^{n_1 + n_2} \times \Gamma$ given by
\[
( a_1,  a_2, \gamma) \efun \gD_1( a_1, \gamma) \gD_2( a_2, \gamma).
\]
It is routine to check that $\gD_1 \otimes \gD_2$ is a predistribution on $\VF^{n_1 + n_2}$ and
\[
\mdl F(\gD_1 \otimes \gD_2) = \wh{\gD}_1 \otimes \wh{\gD}_2.
\]
The \emph{tensor product} of the distributions induced by $\gD_1$ and $\gD_2$ is the distribution induced by $\gD_1 \otimes \gD_2$, which is also written as $\gD_1 \otimes \gD_2$. Note that
\[
\supp(\gD_1 \otimes \gD_2) = \supp(\gD_1) \times \supp(\gD_2).
\]
\end{defn}

If $n_1 = n_2 = n$, we wish to define the convolution $\gD_1 * \gD_2$ of $\gD_1$ and $\gD_2$ as the linear functional on $\mathscr{S}_{n}$ given by
\[
(\gD_1 * \gD_2)(\bm f) = (\gD_1 \otimes \gD_2)(\bm f^{\dag}),
\]
where $\bm f^{\dag}$ is the locally constant integrable function on $\VF^{2n}$ given by $( a_1,  a_2) \efun \bm f( a_1 +  a_2)$. However, as in the classical theory, $\bm f^{\dag}$ is not guaranteed to be a Schwartz-Bruhat function since its support may not be bounded. On the other hand, if either $\gD_1$ or $\gD_2$ has bounded support then $\supp(\bm f^{\dag}) \cap \supp(\gD_1 \otimes \gD_2)$ is bounded. So the classical remedy works:

\begin{defn}
Suppose that $\gD_1$ or $\gD_2$ has bounded support. Then the \emph{convolution} $\gD_1 * \gD_2$ of $\gD_1$ and $\gD_2$ is the linear functional on $\mathscr{S}_{n}$ given by
\[
(\gD_1 * \gD_2)(\bm f) = \int_{( x_1,  x_2) \in \go(0, \gamma)} \bm f( x_1 +  x_2) \bm o^{-2n}_{\beta} (\gD_1 \otimes \gD_2)( x_1,  x_2, \beta),
\]
where $\go(0, \gamma)$ is any definable open polydisc containing $\supp(\bm f^{\dag}) \cap \supp(\gD_1 \otimes \gD_2)$ and $\beta$ is any definable upper bound of $\iota_{\bm f^{\dag} \rest \go(0, \gamma)}(\VF^{2n})$.
\end{defn}

Note that, by Lemma~\ref{pre:dis:vol:par} and Lemma~\ref{dis:supp:bound}, this definition does not depend on the choice of $\gamma$ and $\beta$.

\begin{lem}
The definable function $\VF^n \times \Gamma \fun \KCC$ induced by $\gD_1 * \gD_2$ is a predistribution and hence $\gD_1 * \gD_2$ is a definable distribution.
\end{lem}
\begin{proof}
The proofs of Lemma~\ref{pre:dis:fourier} and Theorem~\ref{dis:fourier} may be easily adapted to work here.
\end{proof}

The following theorem is a generalization of the convolution formula for bounded integrable functions.

\begin{thm}
$\mdl F(\gD_1 * \gD_2) = \wh{\gD}_1 \wh{\gD}_2$.
\end{thm}
\begin{proof}
Let us assume that $\gD_1$ has bounded support. Note that $\wh{\gD}_1 \wh{\gD}_2$ is well-defined since, by Lemma~\ref{dis:bou:four:reg}, $\wh{\gD}_1$ is regular. Choose any locally constant function $\bm h \in \ifn^a_n$ such that $\wh{\gD}_1 = \gD_{\bm h}$. For any $\bm f \in \mathscr{S}_n$, any definable $\gamma \in \Gamma$ that is sufficiently low, and any definable $\beta \in \Gamma$ that is sufficiently large,
\begin{align*}
(\gD_1 * \gD_2)(\wh{\bm f}) &= \int_{( x_1,  x_2) \in \go(0, \gamma)} \Big(\int_{ y \in \go(0, \gamma)} \bm f( y)  \exp_{ x_1 +  x_2}( y) \Big) \bm o^{-2n}_{\beta} \gD_1( x_1, \beta) \gD_2( x_2, \beta)\\
&=\int_{ x_2 \in \go(0, \gamma)} \wh{\gD}_1(\bm f \exp_{ x_2}) \bm o^{-n}_{\beta} \gD_2( x_2, \beta)\\
&= \int_{ x_2 \in \go(0, \gamma)} \Big( \int_{ y \in \go(0, \gamma)} \bm f( y) \exp_{ x_2}( y) \bm h( y) \Big) \bm o^{-n}_{\beta} \gD_2( x_2, \beta).
\end{align*}
Since $\bm f \bm h \rest \go(0, \gamma)$ is a Schwartz-Bruhat function, we deduce that
\[
(\gD_1 * \gD_2)(\wh{\bm f}) = \wh{\gD}_2(\bm f \bm h) =  \wh{\gD}_1 \wh{\gD}_2(\bm f),
\]
as desired.
\end{proof}

\section{The (na\"ive) Weil representation of $\msl_2$}\label{section:weil}

For any self-dual locally compact field $F$, $\cha(F) \neq 2$, and any symplectic vector space $W = V \oplus V^*$ over $F$, where $V$, $V^*$ are two transversal Lagrangian subspaces of $W$, a classical construction due to Shale-Segal-Weil gives a projective representation of the symplectic group $\text{Sp}(W)$ in the Schwartz space of $V$. The goal of this section is to generalize this construction for algebraically closed valued fields of residue characteristic $0$, in a na\"ive way (see the introduction for an explanation), via motivic integration. For simplicity, we shall only do this for $\msl_2$. The canonical reference for such constructions is of course Weil's memoir \cite{weil:rep:64}. On the other hand, Rao's more explicit version \cite{rao:1993} is even more helpful.

\begin{nota}
The notation in Example~\ref{exam:sl2} will be used throughout this section. The action of $n \times n$ matrices on $\VF^n$ is written on the right. For $ x \in \VF^2$ we shall write $\dot x$ for its first coordinate and $\ddot x$ for its second coordinate. The matrix group $\msl_2$ may be presented in a different way. It is generated by elements of the form
\[
u(b) = \begin{pmatrix}
         1 & b\\
         0 & 1
       \end{pmatrix} \in B, \quad
s(a) = \begin{pmatrix}
         a & 0\\
         0 & a^{-1}
       \end{pmatrix} \in A, \quad
w = \begin{pmatrix}
         0  & 1\\
         -1 & 0
       \end{pmatrix} \in W^*,
\]
where $a, b \in \VF$ and $a \neq 0$, subject to the relations
\begin{itemize}
 \item $u$ is an additive homomorphism,
 \item $s$ is a multiplicative homomorphism,
 \item $w s(a) = s(a^{-1}) w$,
 \item $w^2 = s(-1)$,
 \item $w u(a) w = s(- a^{-1}) u(-a) w u(- a^{-1})$.
\end{itemize}
For every $\sigma = \bigl(\begin{smallmatrix}
a & b \\ c & d
\end{smallmatrix}\bigr) \in \msl_2$, if $\sigma \notin B$ then the \emph{standard Bruhat presentation} of $\sigma$ is given by
\begin{align*}
\sigma = u(a/c) w s(-c) u(d/c) = \begin{pmatrix}
         1 & a/c \\
         0 & 1
       \end{pmatrix}\begin{pmatrix}
         0 & 1\\
         -1 & 0
       \end{pmatrix}\begin{pmatrix}
         -c & 0\\
         0 & -1/c
       \end{pmatrix}\begin{pmatrix}
         1 & d / c\\
         0 & 1
       \end{pmatrix}.
\end{align*}
\end{nota}

\begin{defn}
For any quadratic form $q$ on $\VF^n$, the function $\exp_{q} \in \ifn^a_n$ given by $x \efun \exp( \frac{1}{2}q( x) )$ is a \emph{character of $\VF^n$ of second degree}. Let $\rho$ be the $n \times n$ symmetric form associated with $q$ (and hence with $\exp_{q}$). Then
\[
\frac{\exp_{q}(x +  y)}{\exp_{q}( x)\exp_{q}( y)} = \exp(x \cdot  y \rho) = \exp( x \rho \cdot  y).
\]
We say that the character $\exp_{q}$ is \emph{nondegenerate} if $\rho \in \mgl_n$. In that case the function $\exp_{q}' \in \ifn^a_n$ given by $ x \efun \exp_{q}( x \rho^{-1})^{-1}$ is also a character of second degree, which is understood as the dual of $\exp_{q}$ and whose associated symmetric form is clearly $- \rho^{-1}$.

If we write $q'$ for the quadratic form associated with $- \rho^{-1}$ then $\exp_{q}'$ may also be written as $\exp_{q'}$. For any definable subset $A \sub \VF^n$ we write $\exp_{q | A}$ for the restriction $\exp_{q} \rest A$. Note that if $A \sub (\VF^{\times})^n$ is open and bounded then $\exp_{q | A}$ is a Schwartz-Bruhat function. The regular distribution induced by $\exp_{q}$ is written as $\gD_q$.
\end{defn}

Let $\exp_{q}$ be a definable character of $\VF^n$ of second degree, where $n \leq 2$, $q( x) = e x^2$ if $n = 1$, and $q( x) = e_1 \dot x^2 + e_2 \ddot x^2 + e_3 \dot x \ddot x$ if $n=2$. The associated symmetric form is either $\rho = e$ or $\rho = \bigl(\begin{smallmatrix}
e_1 & e_3/2 \\ e_3/2 & e_2
\end{smallmatrix}\bigr)$.

\begin{rem}\label{rem:comp:char2}
Suppose that $n=1$ and $e \neq 0$. Let $b = e /2$ and $\sqrt{b}$ be a chosen square root of $b$. Let $a \in \VF$ and $\gb \sub \VF$ be an open disc around $a$ with $\dim_{\VF}(\gb) = 1$ and $\beta$ the valuative radius of $\gb$. By a computation similar to the one in the proof of Lemma~\ref{simp:char:int}, we have
\[
\int_{x \in \gb} \exp_b(x^2) = \begin{dcases}
               \int_{x \in \OO / \sqrt{b}} \exp_b(x^2), & \text{if } \OO / \sqrt{b} \sub \gb,\\
               \vol(\gb) \exp_b(a^2), & \text{if } \begin{dcases}
                                                \gb \subsetneq \OO / \sqrt{b} \quad \text{or} \\
                                                -\beta - \vv(b) \leq \vv(a) < - \tfrac{\vv(b)}{2},
                                             \end{dcases}\\
               \bm 0, & \text{otherwise};
       \end{dcases}
\]
similarly if $\gb$ is a closed ball. The element
\[
\int_{x \in \OO / \sqrt{b}} \exp_b(x^2) \in \KCC
\]
is denoted by $\bm \theta_b$. If $b =1$ then we simply write $\bm \theta$. Note that if the substructure $S$ contains the square roots of $-1 \in \VF$ then $\bm \theta_b = \bm \theta_{-b}$. The function given by $b \fun \bm \theta_b$ may be considered as a motivic analogue of the Jacobi theta function $\vartheta(0 ; \tau)$ and, in particular, $\bm \theta$ as a motivic analogue of the quadratic Gauss sum. It is conceptually compelling to modify the ring $\KCC$ by various classical identities concerning these objects; but this is not needed in what follows.

Suppose that $n=2$. Let $ a \in \VF^2$ and $\gamma \in \Gamma$. Since $\exp_{q} \rest \go( a, \gamma)$ is an integrable function, if $e_2 \neq 0$ then
\begin{align*}
 \int_{ x \in \go( a, \gamma)} \exp_{q}( x) &= \int_{\dot x \in \go(\dot a, \gamma)} \exp (\tfrac{1}{2}e_1 \dot x^2 ) \int_{\ddot x \in \go(\ddot a, \gamma)} \exp (\tfrac{1}{2}e_2 \ddot x^2 + \tfrac{1}{2}e_3 \dot x \ddot x) \\
 &= \int_{\dot x \in \go(\dot a, \gamma)} \exp  (\tfrac{1}{2} e_1 \dot x^2 ) \int_{\ddot x \in \go(\ddot a + \frac{e_3 \dot x}{2e_2}, \gamma)} \exp \Big(\tfrac{1}{2}e_2 \ddot x^2- \frac{e_3^2\dot x^2}{8 e_2} \Big )\\
 &= \int_{\dot x \in \go(\dot a, \gamma)} \exp \Big(\frac{4 e_1 e_2 - e_3^2}{8 e_2} \dot x^2 \Big) \int_{\ddot x \in \go(\ddot a + \frac{e_3 \dot x}{2 e_2}, \gamma)} \exp_{\frac{e_2}{2}} (\ddot x^2);
\end{align*}
similarly if $e_1 \neq 0$. If $e_1 = e_2 = 0$ then
\[
 \int_{ x \in \go( a, \gamma)} \exp_{q}( x) = \int_{\dot x \in \go(\dot a, \gamma)} \int_{\ddot x \in \go(\ddot a, \gamma)} \exp_{e_3 \dot x}(\ddot x).
\]
Therefore, if $q \neq 0$ and both $\vv(\dot a)$ and $\vv(\ddot a)$ are sufficiently low then, by the third equality above and Lemma~\ref{simp:char:int}, we have $\int_{ x \in \go( a, \gamma)} \exp_{q}( x) = 0$. Of course if $q = 0$ then $\int_{ x \in \go( a, \gamma)} \exp_{q}( x) = \bm o^2_{\gamma}$.
\end{rem}

From now on we suppose that $\exp_q$ is nontrivial, that is, $q \neq 0$, unless indicated otherwise.

\begin{lem}\label{char:conv:b}
For any $\bm f \in \mathscr S_n$, $\bm f * \exp_{q} \in \ifn^b_n$ and hence $\mdl F(\bm f * \exp_{q})$ is well-defined.
\end{lem}
\begin{proof}
The case $n=1$ being completely similar, we shall just show this for the case $n=2$. Since $\bm f$ has bounded support and $\exp_{q}$ is an almost integrable function, an easy argument similar to the one in the proof of Proposition~\ref{convol:C} shows that $\bm f * \exp_{q}$ is also an almost integrable function. So it is enough to show that $\bm f * \exp_{q}$ has bounded support. Let $\go(0, \gamma)$ be a definable open polydisc containing $\supp(\bm f)$ and $\beta > 0$ a definable upper bound of $\iota_{\bm f}(\VF^2)$. For any $ a \in \VF^2$, by the averaging formula,
\[
\int_{ x \in \VF^2} \bm f( x) \exp_{q}( a -  x) = \int_{ x \in \go(0, \gamma)} \Big( \bm o_{\beta}^{-2} \bm f( x) \int_{ y \in \go( x, \beta)} \exp_{q}( a -  y) \Big).
\]
If $\vv(\dot a)$ and $\vv(\ddot a)$ are sufficiently low then, by Remark~\ref{rem:comp:char2},
\[
\int_{ y \in \go( x, \beta)} \exp_{q}( a -  y) = \int_{ z \in \go( a -  x, \beta)} \exp_{q}( z) = \bm 0
\]
and hence $(\bm f * \exp_{q})( a) = \bm 0$.
\end{proof}

The following proposition is an analogue of  \cite[Theorem~2]{weil:rep:64} (also see \cite{cartier:weil:64} and \cite[Theorem~A.1]{rao:1993}), which may be considered as an extension of the convolution formula for integrable functions (see Proposition~\ref{convol:formula:fun}).

\begin{prop}\label{weil:thm:2}
Suppose that $\exp_q$ is nondegenerate. There is a $\bm \gamma_q \in \KCC$, which is called the \emph{Weil index} of $q$, such that, for any $\bm f \in \mathscr S_n$ and any $ a \in \VF^n$,
\[
\mdl F(\bm f * \exp_{q})( a) = \bm \gamma_q \wh{\bm f}( a) \exp_{q'}( a).
\]
\end{prop}
\begin{proof}
The case $n=1$, where $\bm \gamma_q = \bm \theta_{e/2}$, is left to the reader. For $n=2$, we first consider the case $e_1 \neq 0$ or $e_2 \neq 0$. The computations being the same, we just assume that $e_2 \neq 0$. By Lemma~\ref{char:conv:b}, for any $ a$-definable $\gamma < \vv( a)$ that is sufficiently low, we have
\begin{align*}
\MoveEqLeft \mdl F( \bm f * \exp_{q})( a) \\
&= \int_{ x \in \go(0, \gamma)} \Big( \int_{ y \in \go(0, \gamma)} \exp_{q}( y) \bm f( x -  y) \Big) \exp_{ a}( x) \\
&= \int_{ y \in \go(0, \gamma)} \exp_{q}( y) \int_{ z \in \go(0, \gamma)} \bm f( z) \exp_{ a}( z +  y)\\
&= \wh{\bm f}( a) \int_{ y \in \go(0, \gamma)} \exp  (\tfrac{1}{2}e_1 \dot y^2 + \tfrac{1}{2}e_2 \ddot y^2 +  \tfrac{1}{2}e_3 \dot y \ddot y + \dot a \dot y + \ddot a \ddot y  ) \\
&= \wh{\bm f}( a) \int_{\dot y \in \go(0, \gamma)} \exp  (\tfrac{1}{2}e_1 \dot y^2  + \dot a \dot y - \tfrac{e_2}{2}  ( \tfrac{e_3}{2e_2} \dot y + \tfrac{\ddot a}{e_2}  )^2  ) \int_{\ddot y \in \go(\frac{e_3}{2e_2} \dot y + \frac{\ddot a}{e_2}, \gamma)} \exp  (\tfrac{e_2}{2} \ddot y^2 ) .
\end{align*}
By the choice of $\gamma$, we have $\go(\frac{e_3}{2e_2} \dot y + \frac{\ddot a}{e_2}, \gamma) = \go(0, \gamma)$ for all $\dot y \in \go(0, \gamma)$ and hence, by Remark~\ref{rem:comp:char2},
\begin{align*}
\MoveEqLeft \mdl F( \bm f * \exp_{q})( a) \\
&= \wh{\bm f}( a) \bm \theta_{e_2 / 2} \int_{\dot y \in \go(0, \gamma)} \exp_{1/8e_2}((4e_1e_2 - e_3^2) \dot y^2  + 4(2\dot a e_2 - \ddot a e_3) \dot y - 4 \ddot a^2 ) \\
&= \wh{\bm f}( a) \bm \theta_{e_2 / 2} \int_{\dot y \in \go(\frac{2\dot a e_2 - \ddot a e_3}{2 \det(\rho)}, \gamma)} \exp_{\det(\rho)/2 e_2} \Big( \dot y^2  - \frac{\dot a^2 e_2^2 - \dot a \ddot a e_2 e_3 + \ddot a^2 e_1 e_2}{\det(\rho)^2} \Big)
\end{align*}
Again, by the choice of $\gamma$, we have $\go(\frac{2\dot a e_2 - \ddot a e_3}{2 \det(\rho)}, \gamma) = \go(0, \gamma)$ and hence, by Remark~\ref{rem:comp:char2},
\[
\mdl F( \bm f * \exp_{q})( a) = \wh{\bm f}( a) \bm \theta_{e_2 / 2} \bm \theta_{\det(\rho)/2 e_2} \exp_{q'}( a).
\]

Next we consider the case $e_1 = e_2 = 0$.  We have
\begin{align*}
\mdl F( \bm f * \exp_{q})( a) &= \wh{\bm f}( a) \int_{ y \in \go(0, \gamma)} \exp ( \tfrac{1}{2}e_3 \dot y \ddot y + \dot a \dot y + \ddot a \ddot y ) \\
&= \wh{\bm f}( a) \int_{\dot y \in \go(0, \gamma)} \exp (\dot a \dot y) \int_{\ddot y \in \go(0, \gamma)} \exp_{\frac{1}{2} e_3 \dot y + \ddot a}(\ddot y) .
\end{align*}
By the choice of $\gamma$ and Lemma~\ref{simp:char:int}, we deduce
\begin{align*}
\mdl F( \bm f * \exp_{q})( a) &= \wh{\bm f}( a) \bm o_{\gamma} \int_{\dot y \in \gc(- \frac{2 \ddot a}{e_3}, - \gamma - \vv(e_3))} \exp (\dot a \dot y) \\
&=  \wh{\bm f}( a) \bm o_{\gamma} \bm c_{- \gamma - \vv(e_3)} \exp_{q'}( a).
\end{align*}
By Corollary~\ref{dual:vol:cons}, $\bm o_{\gamma} \bm c_{- \gamma - \vv(e_3)} = \bm m_{- \vv(e_3)}$ does not depend on the choice of $\gamma$.

There is an alternative way to compute $\bm \gamma_q$ for the second case, which yields a form that is completely similar to the one in the first case. It goes as follows. By the change of variables $(\dot y, \ddot y) \efun (\dot y, \ddot y - \dot y)$, we have
\begin{align*}
\MoveEqLeft \mdl F( \bm f * \exp_{q})( a) \\
&= \wh{\bm f}( a) \int_{ y \in \go(0, \gamma)} \exp ( \tfrac{1}{2}e_3 \dot y ( \ddot y + \dot y) + \dot a \dot y + \ddot a ( \ddot y + \dot y) ) \\
&= \wh{\bm f}( a) \int_{\ddot y \in \go(0, \gamma)} \exp ( \ddot a \ddot y )  \int_{\dot y \in \go(\frac{1}{2} \ddot y + \frac{\dot a + \ddot a}{e_3}, \gamma)} \exp_{e_3 / 2} \Big(\dot y^2 - \Big(\tfrac{1}{2} \ddot y + \frac{\dot a + \ddot a}{e_3} \Big)^2 \Big)   \\
&= \wh{\bm f}( a) \bm \theta_{e_3 / 2} \int_{\ddot y \in \go(\frac{2(\dot a - \ddot a)}{e_3}, \gamma)} \exp_{- e_3 / 8} \Big( \ddot y^2 +\frac{16 \dot a \ddot a}{e_3^2} \Big)\\
&= \wh{\bm f}( a) \bm \theta_{e_3 / 2} \bm \theta_{\det(\rho)/2 e_3} \exp_{q'}( a),
\end{align*}
where the last two equalities are justified as above.
\end{proof}

\begin{rem}\label{rem:theta}
Since for any $q$ it is easy to arrange $\wh{\bm f}( a) \exp_{q'}( a) = \bm 1$ for some $\bm f$ and some $ a$, we see that, by the computations above, some relations among the elements $\bm \theta_a$ are forced after all: for any $a, b, c \in \VF^{\times}$,
\[
\bm \theta_{a} \bm \theta_{- \frac{a}{4}} = \bm m_{- \vv(a)} \quad \text{and} \quad \bm \theta_{\frac{a}{2}} \bm \theta_{\frac{c}{2a}} = \bm \theta_{\frac{b}{2}} \bm \theta_{\frac{c}{2b}}.
\]
So every $\bm \theta_a$ is a unit. Furthermore, we deduce:
\begin{itemize}
  \item $\bm \theta_{ab} = \bm \theta^{-1} \bm \theta_{a} \bm \theta_{b}$ (in parallel with the classical fact that quadratic Gauss sums are multiplicative) and, in particular, $\bm \theta_{a} = \bm \theta_{- a^2} \bm \theta^{-1} \bm \theta_{- \frac{1}{a}}$ (in parallel with the functional equation $\vartheta(0; a) = (-ia)^{-1/2} \vartheta(0; - 1/a)$ of the Jacobi theta function).
  \item $\bm \theta_a \bm \theta_{-\frac{1}{a}} = \bm 1$ and $\bm \theta_{-a}^2 = \bm \theta \bm \theta_{a^2}$, which imply $\bm \theta^4 = \bm 1$.
  \item $\bm \theta_{\frac{a}{2}}^2 = \bm \theta^2 \bm m_{- \vv(a)}$ and hence $\bm \gamma_q^2 \bm m^2_{\vv(e)/2} = \bm \theta^2$ if $n=1$ (see \cite[Corollary~A.5]{rao:1993} for a comparison).
\end{itemize}
\end{rem}

\begin{rem}\label{rem:2:inv}
Let $\bm f \in \mathscr S_n$ and $\go(0, \gamma)$ be a definable open disc such that $\bm f * \exp_{q | \go(0, \gamma)} = \bm f * \exp_{q}$. Then, by Proposition~\ref{convol:formula:fun},
\[
\mdl F(\bm f * \exp_{q | \go(0, \gamma)}) = \wh{\bm f} \wh{\exp}_{q | \go(0, \gamma)} \in \ifn^b_n.
\]
In particular, if $\bm f = \bm 1_{\gc(0, - \beta)}$ then $\wh{\bm f} = \bm c^n_{- \beta} \bm 1_{\go(0, \beta)}$ and hence, by Proposition~\ref{weil:thm:2}, we see that
\[
\wh{\exp}_{q | \go(0, \gamma)} \rest \go(0, \beta) =  \bm \gamma_q \exp_{q' | \go(0, \beta)}.
\]
\end{rem}

This gives rise to a version of Proposition~\ref{weil:thm:2} in terms of distributions:

\begin{cor}\label{weil:2:cart}
$\wh \gD_{q} = \bm \gamma_q \gD_{q'}$.
\end{cor}
\begin{proof}
It is enough to check that the corresponding predistributions agree. By Theorem~\ref{dis:fourier}, the Plancherel formula (see the paragraph after Theorem~\ref{plancherel}), and Remark~\ref{rem:2:inv}, we have
\[
\wh \gD_{q}(a, \gamma) =  \int_{\VF^n}  \wh{\bm 1}_{\go(0, \gamma)} \exp_q = \bm \gamma_q  \int_{\go(0, \gamma)} \exp_{q'},
\]
as desired.
\end{proof}

\begin{nota}
For each $\sigma = \bigl(\begin{smallmatrix}
a & b \\ c & d
\end{smallmatrix}\bigr) \in \msl_2$ we write $\exp_{\sigma}$ for the character of $\VF^2$ of second degree given by
\[
 x \efun \exp ( \tfrac{1}{2}ab \dot x^2 + \tfrac{1}{2}cd \ddot x^2 + bc \dot x \ddot x ).
\]
The symmetric form associated with $\exp_{\sigma}$ is $\rho = \bigl(\begin{smallmatrix}
ab & bc \\ bc & cd
\end{smallmatrix}\bigr)$. So $\exp_{\sigma}$ is nondegenerate if and only if $\det(\rho) = bc \neq 0$ if and only if $\sigma$ is not in the Borel subgroups $B$, $\overline B$. A simple computation shows that $\exp_{\sigma}'( x) = \exp_{\sigma'} \bigl(\frac{ x}{\sqrt{\det(\rho)}} \bigr)$, where $\sqrt{\det(\rho)}$ is a chosen square root of $\det(\rho)$ and
\[
\sigma' = \begin{pmatrix}
         \frac{-bd}{\sqrt{\det(\rho)}} & \frac{c^2}{\sqrt{\det(\rho)}}\\
         \frac{b^2}{\sqrt{\det(\rho)}} & \frac{-ac}{\sqrt{\det(\rho)}}
       \end{pmatrix} \in \msl_2.
\]
The symmetric form associated with $\exp_{\sigma'}$ is of course $\rho' = - \det(\rho) \rho^{-1}$.

If $\sigma$ is unipotent, that is, if $a = d = 1$ and $bc = 0$, then $\exp_{\sigma}$, written alternatively as $\nu_{b}$ or $\nu_{c}$, will be considered as a character of $\VF$ of second degree. In this case the associated symmetric form $\rho$ is just the scalar $b$ or $c$ and its dual $\rho'$ is just the scalar $-1 / b$ or $-1 /c$.

It is routine to check that $\exp_{\sigma \sigma'} ( x) = \exp_{\sigma} ( x)\exp_{\sigma'} ( x \sigma)$.
\end{nota}

\begin{defn}
For any $ a \in \VF^2$ let $U_{ a}$ be the linear operator on $\mathscr S_1$ given by
\[
U_{ a}(\bm f) : x \efun \exp_{\ddot a}(x) \bm f(x + \dot a).
\]
We clearly have $U_{ a'}U_{ a} = \exp(\ddot a \dot a') U_{ a' +  a}$ and hence the map $ a \efun U_{ a}$ is a projective representation of $\VF^2$.

Let $\sigma = \bigl(\begin{smallmatrix}
a & b \\ c & d
\end{smallmatrix}\bigr) \in \msl_2$. If $\sigma \in B$ then let $R_{\sigma}$ be the linear operator on $\mathscr S_1$ given by
\[
R_{\sigma}(\bm f) : \dot x \efun \bm m_{\vv(a)/2} \exp_{\sigma}(\dot x, \ddot x) \bm f(a\dot x + c \ddot x).
\]
Note that $\ddot x$ does not really occur on the righthand side. We write it in this way for the sake of uniformity, since if $\sigma \notin B$ then we let $R_{\sigma}$ be the linear operator on $\mathscr S_1$ given by
\[
R_{\sigma}(\bm f) : \dot x \efun \int_{\ddot x \in \VF} \exp_{\sigma}(\dot x, \ddot x) \bm f(a\dot x + c \ddot x).
\]
\end{defn}

\begin{lem}\label{weil:end}
Every $R_{\sigma}$ is an endomorphism of $\mathscr S_1$.
\end{lem}
\begin{proof}
This is clear if $\sigma \in B$. Suppose that $\sigma = \bigl(\begin{smallmatrix}
a & b \\ c & d
\end{smallmatrix}\bigr)$ with $c \neq 0$. Let $\bm f \in \mathscr S_1$ and $\go(0, \beta)$ be a definable open disc that contains $\supp(\bm f)$. Then, for all $\dot x \in \VF$, by change of variables and Lemma~\ref{con:vol},
\begin{align*}
R_{\sigma}(\bm f)(\dot x) &= \exp ( \tfrac{1}{2}ab \dot x^2 ) \int_{\ddot x \in \go(- \frac{a}{c}\dot x, \beta - \vv(c))} \exp ( \tfrac{1}{2}cd \ddot x^2 + bc \dot x \ddot x ) \bm f(a\dot x + c \ddot x) \\
&= \bm m_{- \vv(c)} \exp (\tfrac{a}{2c} \dot x^2 ) \int_{\ddot x \in \go(0, \beta)} \exp ( \tfrac{d}{2c} \ddot x^2 - \tfrac{1}{c} \dot x \ddot x ) \bm f(\ddot x).
\end{align*}
If $d = 0$ then $R_{\sigma}(\bm f)(\dot x) = \bm m_{- \vv(c)} \exp (\frac{a}{2c} \dot x^2 ) \wh{\bm f}(- \frac{1}{c} \dot x)$ and hence, by Proposition~\ref{sb:tran:sb}, $R_{\sigma}(\bm f) \in \mathscr S_1$. If $d \neq 0$ then
\[
R_{\sigma}(\bm f)(\dot x) = \bm m_{- \vv(c)} \exp (\tfrac{b}{2d} \dot x^2 ) \int_{\ddot x \in \go(- \frac{1}{d} \dot x, \beta)} \exp_{d/2c} ( \ddot x^2) \bm f (\ddot x + \tfrac{1}{d} \dot x ).
\]
We may choose a definable $\gamma \in \Gamma$ so large that $\bm f (\ddot x + \frac{1}{d} \dot x )$ is constant on all polydiscs of valuative radius $\gamma$. Then, for all open disc $\gb \sub \VF$ of valuative radius $\gamma$ and all $\dot x \in \gb$, the function $\ddot x \efun \bm f (\ddot x + \frac{1}{d} \dot x )$ only depends on $\gb$. This means that $R_{\sigma}(\bm f)$ is locally constant. By the averaging formula, we have
\[
\int_{\ddot x \in \go(- \frac{1}{d} \dot x, \beta)} \exp_{d/2c} ( \ddot x^2) \bm f (\ddot x + \tfrac{1}{d} \dot x ) = \int_{\ddot x \in \go(- \frac{1}{d} \dot x, \beta)} \bm f (\ddot x + \tfrac{1}{d} \dot x ) \bm o^{-1}_{\gamma} \int_{x \in \go(\ddot x, \gamma)} \exp_{d/2c} (x^2)
\]
and hence, by Remark~\ref{rem:comp:char2}, if $\vv(\dot x)$ is sufficiently low then $R_{\sigma}(\bm f)(\dot x) = \bm 0$.
\end{proof}

\begin{lem}\label{weil:cond:1}
For all $ a \in \VF^2$ and all $\sigma \in \msl_2$, $U_{ a} R_{\sigma} = \exp_{\sigma}( a) R_{\sigma} U_{ a \sigma}$.
\end{lem}
\begin{proof}
The computation is straightforward and is left to the reader.
\end{proof}

\begin{lem}\label{weil:cond:2}
For all matrices $\sigma = \bigl(\begin{smallmatrix}
a & b \\ c & d
\end{smallmatrix}\bigr)$,  $\xi_1 = \bigl(\begin{smallmatrix}
a_1 & b_1 \\ 0 & 1/ a_1
\end{smallmatrix}\bigr)$, and $\xi_2 = \bigl(\begin{smallmatrix}
a_2 & b_2 \\ 0 & 1/ a_2
\end{smallmatrix}\bigr)$ in $\msl_2$, we have
\begin{itemize}
 \item $R_{\xi_1}R_{\xi_2} = R_{\xi_1 \xi_2}$,
 \item $R_{\xi_1}R_{\sigma}R_{\xi_2} = \bm m_{\vv(a_1^{-1}a_2)/2} R_{\xi_1 \sigma \xi_2}$ if $c \neq 0$,
 \item $R_{w}R_{\xi_1}R_{w} = \bm \theta_{a_1b_1/2} \bm m_{\vv(a_1b^2_1)/2} R_{w \xi_1 w}$ if $b_1 \neq 0$ and $R_{w}R_{\xi_1}R_{w} = R_{w \xi_1 w}$ if $b_1 = 0$.
\end{itemize}
\end{lem}
\begin{proof}
The computation for the first item is straightforward and the computation for the second item is similar to that for the third item. These are left to the reader. For the third item, let $\bm f \in \mathscr S_1$ and $\go(0, \beta)$ be a sufficiently large definable open disc. If $b_1 \neq 0$ then we have:
\begin{align*}
\MoveEqLeft R_{w}R_{\xi_1}R_{w}(\bm f)(\dot x) \\
&= \int_{\ddot x \in \go(0, \beta)}\exp ( - \dot x \ddot x ) \bm m_{\vv(a_1)/2} \exp (\tfrac{1}{2} a_1 b_1 \ddot x^2 ) \int_{x \in \go(0, \beta)} \exp ( a_1 \ddot x x ) \bm f (- x) \\
&= \bm m_{\vv(a_1)/2} \exp (- \tfrac{1}{2 a_1 b_1}\dot x^2 ) \int_{\ddot x \in \go(0, \beta)} \exp (\tfrac{1}{2} a_1 b_1 \ddot x^2 ) \wh{\bm f} (- a_1 \ddot x - \tfrac{1}{b_1}\dot x )\\
&= \bm \theta_{a_1b_1/2} \bm m_{- \vv(a_1)/2} \exp ( - \tfrac{1}{2 a_1 b_1}\dot x^2 ) \int_{\ddot x \in \go(0, \beta)} \exp (- \tfrac{1}{2a_1 b_1}  \ddot x^2 - \tfrac{1}{a_1 b_1} \dot x \ddot x ) \bm f ( \tfrac{1}{a_1}\ddot x )\\
&= \bm \theta_{a_1b_1/2} \bm m_{\vv(a_1 b^2_1)/2} \int_{\ddot x \in \go(0, \beta)} \exp (- \tfrac{1}{2} a_1 b_1 \ddot x^2 ) \bm f (- \tfrac{1}{a_1} \dot x + b_1 \ddot x )\\
&= \bm \theta_{a_1b_1/2} \bm m_{\vv(a_1b^2_1)/2} R_{w \xi_1 w} (\bm f)(\dot x),
\end{align*}
where the third equality is by the Plancherel formula, the Fourier inversion formula, and Remark~\ref{rem:2:inv}. The case $b_1 = 0$ follows easily from the Fourier inversion formula and is left to the reader.
\end{proof}

\begin{cor}
Let $\sigma = \bigl(\begin{smallmatrix}
a & b \\ c & d
\end{smallmatrix}\bigr)$ and $\sigma' = \bigl(\begin{smallmatrix}
a' & b' \\ c' & d'
\end{smallmatrix}\bigr)$ be matrices in $\msl_2$. The \emph{multiplier} $c(\sigma, \sigma') \in \KCC$ is defined as follows:
\[
c(\sigma, \sigma') = \begin{dcases}
               \bm 1,      &\text{if } c = c' = 0;\\
               \bm m_{- \vv(a)/2},               &\text{if } c = 0 \text{ and } c' \neq 0;\\
               \bm m_{\vv(a')/2},               &\text{if } c \neq 0 \text{ and } c' = 0;\\
               \bm m_{- \vv(cc')/2},               &\text{if } c c' \neq 0 \text{ and } ca' + dc' = 0;\\
               \bm \theta_{\frac{c}{2c'}(ca' + dc')} \bm m_{- \vv(c')} \bm m_{\vv(ca' + dc')},      &\text{otherwise}.
       \end{dcases}
\]
Then we have
\[
R_{\sigma}R_{\sigma'} = c(\sigma, \sigma') R_{\sigma \sigma'}.
\]
Since $c(\sigma, \sigma')$ is always a unit in $\KCC$, it follows that every $R_{\sigma}$ is an automorphism of $\mathscr S_1$. Therefore the map $\sigma \efun R_{\sigma}$ may be considered as \emph{the standard Segal-Shale-Weil (projective) representation} of $\msl_2$.
\end{cor}
\begin{proof}
This follows easily from Lemmas~\ref{weil:end} and \ref{weil:cond:2}. In the computations for the last two cases we can use the standard Bruhat presentations of $\sigma$, $\sigma'$.
\end{proof}

We do not claim that the representation $\sigma \efun R_{\sigma}$ is unique in any sense.

\begin{defn}
Let $\sigma = \bigl(\begin{smallmatrix}
a & b \\ c & d
\end{smallmatrix}\bigr)$ be a matrix in $\msl_2$. The \emph{normalizing constant} $m(\sigma) \in \KCC$ of $\sigma$ is defined as follows:
\[
m(\sigma) = \begin{dcases}
               \bm \theta_{1/2} \bm \theta^{-1}_{d/2} \bm m_{\vv(a)/2},      &\text{if } c = 0;\\
               \bm \theta^{-1}_{c/2},               &\text{otherwise}.
       \end{dcases}
\]
Set $\tilde R_{\sigma} = m(\sigma) R_{\sigma}$ and $\tilde c(\sigma, \sigma') = m(\sigma)m(\sigma') m(\sigma \sigma')^{-1} c(\sigma, \sigma')$.
\end{defn}

\begin{thm}\label{weil:rep}
For all matrices $\sigma = \bigl(\begin{smallmatrix}
a & b \\ c & d
\end{smallmatrix}\bigr)$ and $\sigma' = \bigl(\begin{smallmatrix}
a' & b' \\ c' & d'
\end{smallmatrix}\bigr)$ in $\msl_2$, $\tilde c(\sigma, \sigma') = \bm 1$ and hence the map $\sigma \efun \tilde R_{\sigma}$ is a representation of $\msl_2$.
\end{thm}
\begin{proof}
The computations all rely on the identities in Remark~\ref{rem:theta}. If $c = c' = 0$ then
\[
\tilde c(\sigma, \sigma') = \bm \theta_{d d'/ 2} \bm \theta^{-1}_{d/2} \bm \theta^{-1}_{d'/2} \bm \theta_{1/2} =  \bm 1.
\]
The case $c = 0$ and $c' \neq 0$ is similar. If $c \neq 0$ and $c' = 0$ then
\[
\tilde c(\sigma, \sigma') = \bm \theta_{c a'/ 2} \bm \theta^{-1}_{c/2} \bm \theta^{-1}_{d'/2} \bm \theta_{1/2} \bm m_{\vv(a')} = \bm \theta^{-1}  \bm \theta_{-1} \bm \theta^2_{a'} \bm m_{\vv(a')} =  \bm 1.
\]
If $c c' \neq 0$ and $ca' + dc' = 0$ then
\begin{align*}
\tilde c(\sigma, \sigma') &= \bm \theta_{- c / 2c'} \bm m_{\vv(c / c')/2} \bm \theta_{1/2}^{-1} \bm \theta^{-1}_{c/2} \bm \theta^{-1}_{c'/2}  \bm m_{- \vv(cc')/2}\\
& = \bm \theta^2_{- 1} \bm \theta^2 \bm \theta^{-2}_{1/2} \bm \theta_{1/ c'}^2 \bm m_{- \vv(c')}\\
& = \bm 1.
\end{align*}
Finally, if none of the matrices $\sigma$, $\sigma'$, and $\sigma \sigma'$ is in $B$ then
\begin{align*}
\tilde c(\sigma, \sigma') &= \bm \theta_{\frac{ca' + dc'}{2}} \bm \theta^{-1}_{c/2} \bm \theta^{-1}_{c'/2} \bm \theta_{\frac{c}{2c'}(ca' + dc')} \bm m_{- \vv(c')} \bm m_{\vv(ca' + dc')}\\
& =  \bm \theta^2 \bm \theta_{ca' + dc'}^2 \bm \theta_{1 /c'} \bm \theta_{c'} \bm m_{- \vv(c')} \bm m_{\vv(ca' + dc')}\\
&= \bm \theta \bm \theta_{-1} \bm \theta^4 \\
&= \bm 1.
\end{align*}
Therefore, $\tilde R_{\sigma} \tilde R_{\sigma'} = \tilde R_{\sigma \sigma'}$ in all cases.
\end{proof}

\section{Some technicalities}

\begin{defn}
Let $B$ be a definable subset. A \emph{$\Gamma$-partition} of $B$ is a definable function $\pi : B \fun \Gamma_{\infty}^l$ such that each fiber $B_{\gamma}$ is contained in a (multiplicative) coset of $(\rv^{-1}(\K^{\times}))^n \times (\K^{\times})^m$ and is $\csn(\gamma)$-$\lan{RV}$-definable.
\end{defn}

We say that $B$ is \emph{$\Gamma$-algebraic} if there is a $\Gamma$-partition $\pi$ of $B$ such that each $\pi^{-1}(\gamma)$ is finite.

\begin{lem}\label{gam:alg}
Suppose that there is a definable subset $C \sub B \times \Gamma_{\infty}^n$ such that the fiber $C_{\gamma}$ over any $\gamma \in \Gamma_{\infty}^n$ is finite. Then $B$ is $\Gamma$-algebraic.
\end{lem}
\begin{proof}
By \cite[Lemma~2.21]{Yin:int:expan:acvf} and \omin-minimality, for any $ x \in B$, the fiber $C_{ x}$ contains a definable point $ \gamma_{ x}$. Since the definable function $f : B \fun \Gamma_{\infty}^n$ given by $ x \efun  \gamma_{ x}$ has to be finite-to-one, clearly any $\Gamma$-partition of (the graph of) $f$ witnesses that $B$ is $\Gamma$-algebraic.
\end{proof}

If $B$ is a subset of discs and there is a definable function $f : B \fun \VF$ such that $f(\dot \gb) \in \gb$ for every $\dot \gb \in B$ then $B$ has \emph{definable centers}.

\begin{lem}\label{di:cen}
If $B$ is a $\Gamma$-algebraic subset of discs then it has definable centers.
\end{lem}
\begin{proof}
This is immediate by \cite[Lemma~3.13]{Yin:special:trans} and compactness.
\end{proof}

\begin{defn}\label{defn:di:cen}
Let $A \sub \VF^n  \times \RV^m$. A function $p: A \fun \Gamma$ is an \emph{$\go$-partition} of $A$ if for any $( a,  t) \in A$ the function $p$ is constant on $\go( a, p( a,  t)) \cap A_t$.
\end{defn}

Let $A \sub \VF^n$ and $p : A \fun \Gamma$ be a definable $\go$-partition of $A$. Obviously if $A^* \sub A$ then $p \rest A^*$ is an $\go$-partition of $A^*$ and if $p^*$ is another definable $\go$-partition of $A$ then the definable function on $A$ given by $ a \efun \max \{p( a), p^*( a)\}$ is also an $\go$-partition of $A$. The good behavior of Fourier transform on a Schwartz space, among other things below, will depend on the following lemma:

\begin{lem}\label{vol:par:bounded}
If $A$ is closed and bounded then $p(A)$ is bounded from above.
\end{lem}
\begin{proof}
We first assume that $p$ is $\lan{RV}$-definable. Suppose for contradiction that $p(A)$ is unbounded from above. For each $\gamma \in \Gamma$ let $A_{\gamma} = \set{ a \in A : p( a) > \gamma}$. For each $c \in \MM \mi 0$, by \cite[Corollary~3.14]{Yin:special:trans}, $\acl(c)$ is a model of $\ACVF(S)$ and hence $A_{\vv(c)} \cap \acl(c)$ is nonempty. By compactness, there is a definable subset $B$ of $\MM \times A$ such that $\pr_1(B) = \MM$, $\pr_1 \rest B$ is finite-to-one, and, for each $c \in \MM$ and each $ a \in B_c$, $p( a) > \vv(c)$. Since $\vv(A)$ is bounded from below and $A$ is closed, by \cite[Lemma~9.5]{Yin:special:trans}, there is a finite subset $L \sub A$ such that $\lim_{\MM \mi \set{0} \rightarrow 0} B = L$ (for the meaning of this notation see \cite[Definition~9.1]{Yin:special:trans}). So, for any $ a \in L$, there is a $c \in \MM$ with $\vv(c) > p( a)$ such that $B_c \cap \go( a, p( a)) \neq \0$. Since $p(B_c) >  \vv(c)$, this contradicts the assumption that $p$ is an $\go$-partition of $A$.

We proceed to the general case. Note that, since $p$ is an $\go$-partition of $A$ and $A$ is bounded, it follows that $p(A)$ is also bounded from below, that is, (the graph of) $p$ is bounded. Also, by \cite[Corollary~2.23]{Yin:int:expan:acvf}, $I = \vv(A)$ is $\lan{RV}$-definable. Let $\pi : p \fun \Gamma_{\infty}^l$ be a $\Gamma$-partition of $p$. It is easy to see that, since $p$ is bounded, we may assume that $\pi$ is bounded as well. By \cite[Corollary~2.23]{Yin:int:expan:acvf}, quantifier elimination, and compactness, we may also assume that there are an $\lan{RV}$-definable finite partition $D_i$ of $\pi(p)$ and finitely many terms $\tau_i( X,  Y)$ of the form
\[
\vv(F_i(X,  Y)) - \vv(G_i( X, Y)) + \alpha_i,
\]
where $F_i( X,  Y)$, $G_i( X,  Y)$ are polynomials and $\alpha_i \in \Gamma(S)$, such that, for each $ \gamma \in D_i$, the function $\pi^{-1}( \gamma)$ is given by $ a \efun \tau_i( a, \csn( \gamma))$. Let $B_i \sub \vv^{-1}(I \times D_i) \times \Gamma$ be the subset defined by the formula
\[
Z= \tau_i( X,  Y) \wedge F_i( X,  Y) \neq 0 \wedge G_i( X,  Y) \neq 0.
\]
Set $B = \bigcup_i B_i$, which is the graph of an $\lan{RV}$-definable function $p' : A' \fun \Gamma$ such that $\pi^{-1}( \gamma) \sub p' \rest \fib(A', \csn( \gamma))$ for every $ \gamma \in \pi(p)$. Observe that if $p'( a,  b) = \gamma$ then there is a $\delta_{ a,  b} \in \Gamma$ such that
\[
\go(( a,  b), \delta_{ a,  b}) \sub A' \quad \text{and} \quad p'(\go(( a,  b), \delta_{ a,  b})) = \{\gamma\}.
\]
By \omin-minimality, we may assume that $\delta_{ a,  b}$ is $( a,  b)$-$\lan{RV}$-definable and is the least element that makes this hold. Then the function $p'' : A' \fun \Gamma$ given by $( a,  b) \efun \max \{p'( a,  b), \delta_{ a,  b}\}$ is an $\lan{RV}$-definable $\go$-partition of $A'$.

For each $ \gamma \in \pi(p)$, the topological closure $A_{ \gamma}$ of $\dom(\pi^{-1}( \gamma))$ is a $\csn( \gamma)$-$\lan{RV}$-definable subset of $A$. So there is a $\csn( \gamma)$-$\lan{RV}$-definable $\epsilon_{ \gamma} \in \Gamma$ such that $p(A_{ \gamma}) < \epsilon_{ \gamma}$. By \cite[Corollary~2.23]{Yin:int:expan:acvf} and compactness, the function $h : \pi(p) \fun \Gamma$ given by $ \gamma \efun \epsilon_{ \gamma}$ is $\lan{RV}$-definable. Now let
\[
 A'' = A' \mi \bigcup_{( a,  b) \in \vv^{-1}(I \times \pi(p)) \mi A'} \go(( a,  b), h(\vv( b))).
\]
Clearly $A''$ is $\lan{RV}$-definable and closed. Moreover, by the construction of $h$, we still have $\pi^{-1}( \gamma) \sub p' \rest \fib(A'', \csn( \gamma))$ for every $ \gamma \in \pi(p)$. Since $p''(A'')$ is bounded from above, we see that $p(A)$ is bounded from above as well.
\end{proof}

Recall \cite[Notation~3.16, Definition~8.8]{Yin:special:trans} and \cite[Definitions~2.16, 4.18]{Yin:int:expan:acvf}. In $\gC^{\csn}$, \cite[Lemma~8.9]{Yin:special:trans} still holds:

\begin{lem}\label{fun:loc:cons}
Every definable function $f : \VF^n \fun \mdl P(\RV^m)$ is locally constant almost everywhere.
\end{lem}
\begin{proof}
Let $\phi( X,  Y)$ be a quantifier-free formula that defines $f$. Let $\tau_i( X,  Y)$ enumerate the top occurring $\VF$-terms of $\phi$. By \cite[Theorem~4.25]{Yin:int:expan:acvf} and compactness, there is a special bijection $T$ on $\VF^n$ such that, for every $\rv$-polydisc $\gp \sub T(\VF^n)$, every $ t \in \RV^m$, and every $i$, the function
\[
\rv \circ \tau_i(-,  t) \circ (T^{-1} \rest \gp)
\]
is constant. This means that, for all $ a,  b \in T^{-1} (\gp)$, $f( a) = f( b)$. If $\dim_{\VF}(\gp) = n$ then, by \cite[Lemma~4.7]{Yin:int:expan:acvf} and \cite[Lemma~4.6]{Yin:special:trans}, $T^{-1}(\gp)$ contains an open polydisc. The lemma follows.
\end{proof}

\begin{defn}\label{defn:dual}
We say that two polydiscs $\gb, \gd \sub \VF^n$ have the same \emph{signature} if, for all $i \leq n$, $\dim_{\VF}(\pr_i(\gb)) = \dim_{\VF}(\pr_i(\gd))$ and $\pr_i(\gb)$ is an open disc if and only if $\pr_i(\gd)$ is an open disc.

For any $\eta \in \Gamma$, the \emph{$\eta$-dual disc} of $\go(a, \gamma)$ is $\gc(a, - \gamma + \eta)$ and the \emph{$\eta$-dual disc} of $\gc(a, \gamma)$ is $\go(a, - \gamma + \eta)$. For any $ \eta \in \Gamma^n$, the \emph{$ \eta$-dual polydisc} of a polydisc $\gb \sub \VF^n$ around $(a_1, \ldots, a_n)$ is the polydisc $\gb^{\iota( \eta)}$ around $(a_1, \ldots, a_n)$ such that each $\pr_i (\gb^{\iota( \eta)})$ is the $\eta_i$-dual of $\pr_i (\gb)$ with respect to $a_i$. If $ \eta = 0$ then it will be dropped from the notation.
\end{defn}

\begin{cor}\label{dual:vol:cons}
Suppose that $ \eta \in \Gamma^n$ is definable and the definable polydiscs $\gb, \gd \sub \VF^n$ have the same signature. Then
\[
\vol(\gb)\vol(\gb^{\iota( \eta)}) = \vol(\gd)\vol(\gd^{\iota( \eta)}).
\]
In particular, $\vol(\gb)\vol(\gb^{\iota}) =  \bm 1$ and hence both $\vol(\gb)$ and $\vol(\gb^{\iota})$ are units in $\KRC$.
\end{cor}
\begin{proof}
Clearly it is enough to show the case $n = 1$. By Lemma~\ref{di:cen}, $\gb$, $\gd$ have definable centers. So, without loss of generality, we may assume that they are both centered around $0$. Let $\beta$, $\delta$ be the valuative radii of $\gb$, $\gd$. Now, using the matrix $\biggl(\begin{smallmatrix}
\frac{\csn(\beta)}{\csn(\delta)} & 0 \\ 0 & \frac{\csn(\delta)}{\csn(\beta)}
\end{smallmatrix}\biggr)$, the claim follows immediately from the change of variables formula.
\end{proof}

\begin{nota}\label{nota:vol:m}
For any definable $\gamma, \eta \in \Gamma$, since $\bm o_{\gamma + \eta} / \bm o_{\gamma} = \bm o_{\gamma + \eta} \bm c_{-\gamma}$ does not depend on $\gamma$, we may denote it by $\bm m_{\eta}$. If $\eta' \in \Gamma$ is also definable then
\[
\bm o_{\gamma + \eta + \eta'} \bm o_{\gamma} = \bm o_{\gamma + \eta} \bm o_{\gamma + \eta'}
\]
and hence $\bm m_{\eta + \eta'} = \bm m_{\eta}\bm m_{\eta'}$. This implies that $\bm m_{\eta}$ is invertible in $\KRC$ and $\bm m_{\eta}^{-1} = \bm m_{- \eta}$. We also have $\bm o_{\gamma + \eta} \bm c_{\gamma} = \bm o_{\gamma} \bm c_{\gamma + \eta}$ and hence $\bm m_{\eta} = \bm c_{\gamma + \eta} / \bm c_{\gamma}$.
\end{nota}

%

\begin{thebibliography}{10}

\bibitem{bernstein:joseph:1972}
Joseph Bernstein, \emph{Analytic continuation of distributions with respect to
  a parameter}, Functional Analysis and its Applications \textbf{6} (1972),
  no.~4, 26--40.

\bibitem{cartier:weil:64}
Pierre Cartier, \emph{{\"U}ber einige integralformeln in der theorie der
  quadratischen formen}, Mathematische Zeitschrift \textbf{84} (1964), no.~2,
  93--100.

\bibitem{cluckers:hales:loeser:transfer}
Raf Cluckers, Thomas Hales, and Fran\c{c}ois Loeser, \emph{Transfer principle
  for the fundamental lemma}, On the stabilization of the trace formula
  (L.~Clozel, M.~Harris, J.-P. Labesse, and B.-C. Ng{\^o}, eds.), International
  Press of Boston, 2011, arXiv:0712.0708v1.

\bibitem{cluckers:loeser:constructible:motivic:functions}
Raf Cluckers and Fran\c{c}ois Loeser, \emph{Constructible motivic functions and
  motivic integration}, Inventiones Mathematicae \textbf{173} (2008), no.~1,
  23--121, math.AG/0410203.

\bibitem{cluckers:loeser:motivic:fourier:transform}
\bysame, \emph{Constructible exponential functions, motivic {F}ourier transform
  and transfer principle}, Annals of Mathematics \textbf{171} (2010), no.~2,
  1011--1065, math.AG/0512022.

\bibitem{cunning:hales:good}
Clifton Cunningham and Thomas Hales, \emph{Good orbital integrals},
  Representation Theory \textbf{8} (2004), 414--457, arXiv:math/0311353v2.

\bibitem{gelfand:shilov:1964}
I.~M. Gel'fand and G.~E. Shilov, \emph{Generalized functions: properties and
  operations}, vol.~1, Academic Press, New York, 1964, translated from the
  Russian by Eugene Saletan.

\bibitem{gordon:yaffe:2008}
Julia Gordon and Yoav Yaffe, \emph{An overview of arithmetic motivic
  integration}, Ottawa lectures on admissible representations of reductive
  {\p-adic} groups, Fields Institute Monograph, vol.~26, Amer. Math. Soc.,
  Providence, RI, 2009, arXiv:0811.2160v1.

\bibitem{iwo:alg:HPK}
Thomas~J. Haines, Robert~E. Kottwitz, and Amritanshu Prasad,
  \emph{{I}wahori-{H}ecke algebras}, J. Ramanujan Math. Soc. \textbf{25}
  (2010), no.~2, 113--145, arXiv:math/0309168v3.

\bibitem{hales:2005}
Thomas Hales, \emph{What is motivic measure?}, Bulletin of the American
  Mathematical Society \textbf{42} (2005), no.~2, 119--135, math.AG/0312229.

\bibitem{hormander:83}
Lars H\"{o}rmander, \emph{The analysis of linear partial differential operators
  {I}: distribution theory and {F}ourier analysis}, Grundlehren der
  mathematischen Wissenschaften 256, Springer-Verlag, New York, 1983.

\bibitem{hrushovski:kazhdan:integration:vf}
Ehud Hrushovski and David Kazhdan, \emph{Integration in valued fields},
  Algebraic geometry and number theory, Progr. Math., vol. 253, Birkh\"{a}user,
  Boston, MA, 2006, math.AG/0510133, pp.~261--405.

\bibitem{hru:kazh:val:ring:2006}
\bysame, \emph{The value ring of geometric motivic integration and the
  {I}wahori {H}ecke algebra of $\textup{SL}_2$}, Geometric and Functional
  Analysis GAFA \textbf{17} (2008), no.~6, 1924--1967, with an appendix by Nir
  Avni, math.AG/0609115.

\bibitem{hru:kazh:2009}
\bysame, \emph{Motivic {P}oisson summation}, Mosc. Math.~J. \textbf{9} (2009),
  no.~3, 569--623, arXiv:0902.0845v1.

\bibitem{loo:mot}
Eduard Looijenga, \emph{Motivic measures}, {S\'eminaire Bourbaki, Volume
  1999/2000, Expos\'es 865-879, Paris: Soci\'et\'e Math\'ematique de France,
  Ast\'erisque 276, 267-297, Exp. No. 874 (2002), arXiv:math/0006220v2}.

\bibitem{rao:1993}
R.~Ranga Rao, \emph{On some explicit formulas in the theory of {W}eil
  representation}, Pacific Journal of Mathematics \textbf{157} (1993), no.~2,
  335--372.

\bibitem{weil:rep:64}
Andr\'{e} Weil, \emph{Sur certains groupes d'op{\'e}rateurs unitaires}, Acta
  Mathematica \textbf{111} (1964), no.~1, 143--211.

\bibitem{weil:basic:number}
\bysame, \emph{Basic number theory}, 3rd ed., Springer-Verlag, New York, 1974.

\bibitem{Yin:QE:ACVF:min}
Yimu Yin, \emph{Quantifier elimination and minimality conditions in
  algebraically closed valued fields}, arXiv:1006.1393v1, 2009.

\bibitem{Yin:special:trans}
\bysame, \emph{Special transformations in algebraically closed valued fields},
  Annals of Pure and Applied Logic \textbf{161} (2010), no.~12, 1541--1564,
  arXiv:1006.2467.

\bibitem{Yin:int:acvf}
\bysame, \emph{Integration in algebraically closed valued fields}, Annals of
  Pure and Applied Logic \textbf{162} (2011), no.~5, 384--408,
  arXiv:0809.0473v2.

\bibitem{Yin:int:tcvf}
\bysame, \emph{Additive invariants in \omin-minimal valued fields},  (2013),
  arXiv:1307.0224, submitted.

\bibitem{Yin:int:expan:acvf}
\bysame, \emph{Integration in algebraically closed valued fields with
  sections}, Annals of Pure and Applied Logic \textbf{164} (2013), no.~1,
  1--29, arXiv:1204.5979v2.

\end{thebibliography}

\end{document}